\documentclass[11pt]{article}

\usepackage{natbib} 
\usepackage{color} 
\usepackage{hyperref} 
\hypersetup{
    colorlinks=true,
    linkcolor=red,
    citecolor=blue,
    urlcolor=blue,
    pdfborder={0 0 0}
}
\usepackage{authblk}        
\usepackage[utf8]{inputenc} 
\usepackage[T1]{fontenc}    
\usepackage{url}            
\usepackage{booktabs}       
\usepackage{amsfonts}       
\usepackage{lipsum}         
\usepackage{geometry}       
\usepackage{algorithm}
\usepackage{algpseudocode}

\usepackage{amsmath, amssymb, amsthm} 
\usepackage{rsfso}      
\usepackage[cal=cm,        
scr=rsfso,                   
calscaled=.96                
]{mathalfa}                  
\usepackage{dsfont}
\usepackage{comment}     

\usepackage{tikz}

\usepackage[toc,page]{appendix}  

\newtheorem{proposition} {Proposition} [section]
\theoremstyle{definition}
\newtheorem{definition}{Definition}[section]
\newtheorem{example}{Example}[section]
\theoremstyle{remark}
\newtheorem{remark}{Remark}[section]

\newcommand{\argmin}{\operatornamewithlimits{argmin}}

\newcommand{\bs}[1]{\boldsymbol{#1}}
\newcommand{\mds}[1]{\mathds{#1}}

\newcommand{\fr}[2][n]{\frac{{#2}}{{#1}}}
\newcommand{\parenf}[1]{\left(#1\right)}
\newcommand{\parenc}[1]{\left\{#1\right\}}
\newcommand{\parens}[1]{\left[#1\right]}

\newcommand{\frsum}{\fr{1}\sum_{i=1}^n}
\newcommand{\tr}{\text{tr}}

\providecommand{\keywords}[1]{\textbf{\textit{Keywords:}} #1}

\title{Robust Estimation of Sparse, High Dimensional Time Series with Polynomial Tails.}

\author[1]{Sagnik Halder\thanks{Email:shalder@ufl.edu}}
\author[2]{George Michailidis\thanks{Email: gmichail@ufl.edu}}
\affil[1]{Department of Statistics, University of Florida}
\affil[2]{Department of Statistics \& Informatics Institute, University of hFlorida}

\usepackage{xparse}

\ExplSyntaxOn

\NewDocumentCommand{\definealphabet}{mmmm}
 {
  \int_step_inline:nnn { `#3 } { `#4 }
   {
    \cs_new_protected:cpx { #1 \char_generate:nn { ##1 }{ 11 } }
     {
      \exp_not:N #2 { \char_generate:nn { ##1 } { 11 } }
     }
   }
 }

\ExplSyntaxOff

\definealphabet{mb}{\mathbb}{A}{Z}
\definealphabet{mc}{\mathcal}{A}{Z}
\definealphabet{mr}{\mathrm}{A}{Z} 
\definealphabet{mf}{\mathfrak}{A}{Z}
\definealphabet{mf}{\mathfrak}{a}{z}
\definealphabet{ms}{\mathscr}{A}{Z} 
\definealphabet{msf}{\mathsf}{A}{Z} 
\definealphabet{msf}{\mathsf}{a}{z}

\begin{document}

\maketitle

\begin{abstract}
High dimensional Vector Autoregressions (VAR) have received a lot of interest recently due to novel applications in health, engineering, finance and the social sciences. Three issues arise when analyzing VAR$s$: (a) The high dimensional nature of the model in the presence of many time series that poses challenges for consistent estimation of its parameters;  (b) the presence of \textit{temporal} dependence introduces additional challenges for theoretical analysis of various estimation procedures; (b) the presence of heavy tails in a number of applications. Recent work, e.g. \citep{basu2015regularized},\citep{kock2015oracle}, has addressed consistent estimation of \textit{sparse} high dimensional $\textit{stable}$ $\textit{Gaussian}$ VAR models based on an $\ell_1$ LASSO procedure. Further, the rates obtained are optimal, in the sense that they match those for iid data, plus a multiplicative factor (which is the "price" paid) for temporal dependence.

However, the third issue remains unaddressed in extant literature. 
This paper extends existing results in the following  important direction: it considers consistent estimation of the parameters of sparse high dimensional VAR models driven by heavy tailed  homoskedastic or heteroskedastic noise processes (that do not possess all moments). 
A \textit{robust} penalized approach (e.g., LASSO) is adopted 
for which optimal consistency rates and corresponding \textit{finite} sample bounds for the underlying model parameters are obtain that \textit{match} those for iid data, albeit paying a price for temporal dependence. The theoretical results are illustrated on VAR models and also on other popular time series models. Notably, the key technical tool used, is a single $\textit{concentration bound}$ for heavy tailed dependent processes.
\end{abstract}


\section{Introduction}\label{intro}

A number of statistical models are routinely used in modeling and analysis of temporally dependent high dimensional data encountered in neuroscience, genomics, economics, finance and signal processing applications. To obtain consistent estimates of the underlying model parameters and thus overcome the limited availability of samples (time points) --- the true, underlying parameter space is assumed to be low-dimensional or "sparse". This sparsity, in turn, is recovered by penalized ("regularized") methods, with popular methods being the LASSO or the  "Dantzig" selector \citep{wainwright2019high}. To set the stage, consider a stochastic regression model \citep{hamilton2020time} given by
\begin{gather}\label{stoch-reg}
y_t=x_t^\top\beta^*+\epsilon_t,\quad t=1,\dots,n,    
\end{gather}
where $\{y_t\}$ and $\{x_t\}$ are generated from stationary stochastic processes in $\mbR$ and $\mbR^p$, respectively, and $\{\epsilon_t\}$ is the noise process. Writing $Y^\top=[y_1\dots y_n],X^\top=[x_1 \dots x_n]$ and $\mcE^\top=[\epsilon_{1} \dots \epsilon_{n}]$, the model can be expressed in matrix form as 
\begin{gather} \label{stoch-reg-matrix-form}
Y = X \beta^*+ \mcE .
\end{gather}
The goal is to consistently estimate the $p$-dimensional parameter $\beta^*$ from data pairs $(y_t,x_t)$, where  $1\le t\le n$, when the dimension $p\gg n$. To this end, the standard assumption is that $\beta^*$ has an underlying low-dimensional structure, e.g. being sparse or group sparse \citep{wainwright2019high}, which can be induced using a penalized estimation approach.

The main application of this basic setup will be a VAR model of lag $d$ in $\mbR^p$. Formally, a VAR($d$) process ($d\ge 1$ represents a $\textit{fixed}$ lag) is defined as follows: $Z_t=B_1^\top Z_{t-1}+\dots +B_d^\top Z_{t-d} +\epsilon_t$, where each $B_k, k=1,\dots,d$ is a fixed coefficient matrix in $\mbR^{p\times p}$ and innovations $\{\epsilon_t\}$ are $p$-dimensional random vectors whose components $\{\epsilon_{tj}\}$ are independent random variables. The parameter of interest is the matrix $B^\top=[B_1^\top,\dots,B_d^\top]$ with $dp^2$ entries. It can be seen that even for small values of $p=20$, $d=5$, there are 2000 parameters and in many applications the available number of time points is significantly smaller. Note that every VAR(d) process has an equivalent VAR(1) representation \cite{lutkepohl2005new} as $\tilde{Z}_t =\tilde{B}^\top\tilde{Z}_{t-1} +\tilde{\epsilon}_t$ where
\begin{gather}\label{eq:VAR1}
\tilde{Z_t}:=
\begin{bmatrix}
Z_t \\ 
Z_{t-1} \\
\vdots \\
Z_{t-d+1}
\end{bmatrix}_{(pd \times 1)} ,
\quad
\tilde{\epsilon}_t:=
\begin{bmatrix}
\epsilon_t \\ 
\bs{0} \\
\vdots \\
\bs{0} 
\end{bmatrix}_{(pd \times 1)}  ,
\quad
\tilde{B}^\top:=
\begin{bmatrix}
B_1^\top & B_2^\top & \cdots & B_{d-1}^\top & B_d^\top \\ 
I_p & \bs{0} & \cdots & \bs{0}  & \bs{0} \\ 
\bs{0} & I_p &        & \bs{0}  & \bs{0} \\ 
\vdots &     & \ddots &  \vdots & \vdots\\ 
\bs{0} & \bs{0} & \cdots & I_p & \bs{0} 
\end{bmatrix}_{(dp\times dp)} .
\end{gather}
Suppose that observations for $T\ge d$ data points $\{Z_0,\dots,Z_T\}$ are available; we can then write the original VAR(d) process as $Z_t=B^\top\tilde{Z}_{t-1}+\epsilon_t,d\le t\le T$. Splitting this up into $p$ parallel regressions, the $j^{th}$ regression is given by
\begin{gather}
Z_{tj}= B_{j:}^\top\tilde{Z}_{t-1}+\epsilon_{tj}, \quad d\le t\le T,    
\end{gather}
wherein $ B_{j:}$ represents the $j^{th}$ row of $B^\top$. Fix $1\le j\le p$ and rewrite $Z_{tj}=y_{t-d+1}$, $B_{j:}=\beta^*$, $\tilde{Z}_{t-1}=x_{t-d+1}$ and $\epsilon_{tj}=\eta_{t-d+1}$ (we suppress dependence on $j$ for now). Hence, the $j^{th}$ regression is just a stochastic regression $y_i=x_i^\top\beta^* +\eta_i,1\le i\le n$, where $n=T-d+1$, or 
\begin{gather}\label{VAR-stoch-reg}
Y=X\beta^*+\eta    
\end{gather}
in matrix form, and the parameter $\beta^*\in\mbR^{pd}$.
Clearly, this is an example of the generic stochastic regression framework \eqref{stoch-reg}. Hence, we can carry out $p$ penalized regressions $\textit{simultaneously}$ with the $\textit{same}$ penalty $\mcR$ and some tuning parameter $\lambda_n$ (as given by \eqref{robust-loss}). The only difference will be the error probability (regret bound) for simultaneous estimation: if the regret bound for each regression is $\mfp_{err}$ (which will be the same for every regression as it only depends on $\mcR$), then the consistency result holds with probability at least $p(1-\mfp_{err})-(p-1)=1- p\mfp_{err}$ (simple Bonferroni bound). 

One key feature of this model is that, it is an endogenous process---the predictors $\{x_t\}$ and output $\{y_t\}$ are both driven by the noise $\{\eta_t\}$. As a result two things happen: (a) if the noise $\{\eta_t\}$ is non-Gaussian/heavy tailed, this causes both $\{x_t\}$ and $\{y_t\}$ to be non-Gaussian/heavy tailed, and (b) the temporal dependence in the predictors $\{x_t\}$ is due to noise $\{\eta_t\}$ and further, is a function of the regression parameter $\beta^*$.

Two popular penalized methods have been extensively studied to date --- the $\textit{Dantzig}$ estimator and $\textit{LASSO}$. Further, in order to induce a general sparsity pattern, penalized estimation using different choices of penalty (usually a norm), have been extensively studied. However, their properties have been usually studied under the fixed design setup ($\{x_t\}$ is deterministic) or the iid case, where the $\{x_t\}$ typically has iid Gaussian or sub-Gaussian entries/rows, but correlated columns to signify cross-sectional dependence between the $p$ components of each $\{x_t\}$. 

In time series data, a key challenge is to handle $\textit{temporal}$ dependence, in addition to cross sectional dependence. Hence, the methods used in extant literature are not directly applicable in our case and require a natural extension. For a comprehensive overview on recent contributions in high dimensional time series, see \citep*[section 1.1]{wong2020lasso}. The main difficulty, specifically, is to apply concentration bounds for sums of dependent random variables, and get (a) tight bounds, i.e, those that match the iid bounds, e.g., those of Hoeffding or Bernstein inequalities, and (b) interpretable or explicit constants that clearly display the temporal dependence. In our case, such a concentration bound, is given by Proposition \eqref{ConcMarkov}. This is a concentration bound for bounded functions of sums of random variables generated from a stationary Markov Chain. Note that there exists a substantial body of literature on concentration bounds of this type. However, the majority of the results are devoted to (a) \textit{finite} state space, (b) \textit{uniformly} mixing/ergodic or (c) \textit{reversible} Markov chains \citep{meyn2012markov}. However, none of those apply to the heavy tailed VAR model under consideration --- e.g. see \citep[Theorem 2]{athreya1986mixing} for uniformly mixing autoregressions, and \citep[Section 3.5]{adamczak2008tail} or \citep[Remark 4]{adamczak2015exponential}. Moreover, when those conditions are applicable to other time series models of interest, they do not usually give rise to easily interpretable constants, or are not tight \citep[Remark 2]{adamczak2015exponential}. Therein, lies the main difficulty and challenge that the current paper addresses.

An additional complication arises when the noise $\{\epsilon_t\}$ is heavy tailed (having only finitely many moments), which in turn, induces noisiness in $\{x_t\}$ and $\{y_t\}$ as well (e.g. in autoregressions). In this case, a robust regression framework with a general penalty, is indispensable. A gold standard in classical robust statistics is Huber regression (or M-estimation in general), which has recently been extended to the high-dimensional setting \citep{loh2017statistical, loh2018scale}. We leverage this framework for the stochastic regression model \eqref{stoch-reg}, in order to consistently estimate the regression parameter $\beta^*$ at optimal rates. Finally, we apply our results to a large class of high dimensional, heavy tailed, VAR models (with or without $\textit{heteroskedastic}$ noise), under general sparsity, for which few theoretical results exist thus far. The main contributions of this work are summarized next:

\begin{enumerate}
    \item We extend theoretical consistency results in high-dimensional literature from the iid setup to the case of temporal dependence under a general sparsity pattern.
    \item Provide non-asymptotic consistency rates and finite sample bounds for a robust approach to heavy tailed, temporally processes. These rates match the optimal rates obtained in the iid case, in many time series models.
    \item Apply our theoretical results to several interesting examples of both linear and non-linear high dimensional time series.
\end{enumerate}

\subsection{Notation}

Throughout the paper, we use the following notation: $||\cdot||$ denotes the $\ell_2$-norm of a vector, while $||\cdot||_p$, $||\cdot||_F$ and $||\cdot||_{nuc}$ denote the matrix norm induced by the $\ell_p$ norm $(1\le p\le\infty)$, the Frobenius  and nuclear norms of a matrix, respectively. The sparsity (or number of non-zero entries) of a matrix is denoted by $||\cdot||_0$, while$||\cdot||_{max}$ denotes the maximum of absolute values of entries of a matrix. For a $p\times d$ matrix A and $G\subseteq \{1,...,p\}\times \{1,...,d\}$ we write the submatrix $A_G = [A_{i,j} : (i,j) \in G]$. For a partition of the set $\{1,...,p\}\times \{1,...,d\}$ into disjoint groups $G_1,...,G_M$, we denote the group norm of a matrix $A$ as $||A||_{2,1} = \sum_{i=1}^M ||A_{G_i}||_F$. We also denote the mixed norm $||A||_{2,\infty}=\underset{1\le j\le d}{max} ||A_{.j}||$, where $A_{i.}$ and $A_{.j}$ denote the $i^{th}$ row and $j^{th}$ column of $A$, respectively. Let $\mbB_\mcR (0,1)=\{v:\mcR(v)\le 1\}$ and write $\mbB_2$ for $\mbB_\mcR (0,1)$ when $\mcR=||\cdot||$.  The dual of norm $\mcR()$ will be denoted by $\mcR^*()$. The cardinality of a set $J$ is denoted by $card(J)$, while its closure and convex hull is denoted by $cl(J)$ and $conv(J)$, respectively. We use $\{e_1 , \dots e_p \}$ to denote the standard canonical vectors in $\mbR^p$ with respect to the $\ell_2$ norm. For positive real numbers $a,b$, we write $b \succsim a$ if there is a positive constant $c$, independent of the model parameters, such that $b\ge ca$. Also, we say $a\asymp b$ if $a \succsim b$ and $b \succsim a$. Absolute model-free positive constants are usually denoted by $c_i$ and may change from line to line throughout the paper (these are of little concern, since they do not impact the results). The conjugate of a complex matrix $A$ is denoted as $A^*$, and if $A$ is a real matrix, its transpose is written as $A^\top$. The maximum and minimum eigenvalues of a matrix $A$ are written as $\Lambda_{\max}(A)$ and $\Lambda_{\min}(A)$. The trace and determinant of a square matrix $A$ is $\tr(A)$ and $\det(A)$, respectively.

\section{Setting Up the Problem: Preliminaries.}
 
We apply the following penalized estimation procedure (a variant of LASSO) to the model \eqref{stoch-reg}:
\begin{gather}\label{robust-loss}
\hat{\beta}\in\underset{\beta\in\mbR^p}{\argmin} \parenc{\mcL_n(\beta) +\lambda_n\mcR(\beta)},   
\end{gather}
wherein $\mcL_n(\beta)=\frsum\mfw(x_i)\ell[\mfw(x_i) (y_i- x_i^\top\beta)]$ denotes a "robust" objective function, and $\mfw()$ is a weight function which dampens the noisiness of the predictor observations $x_1,\dots,x_n$. In case $\ell(\cdot)$ is the least squares function, and the weight function $\mfw\equiv 1$, the procedure \eqref{robust-loss} is just the usual LASSO. Here $\mcR()$ is the penalty function (usually a norm) and $\lambda_n$ is the tuning parameter. Before proceeding further, we need to unpack some notions on (a) the function $\ell()$, (b) the penalty $\mcR()$ and (c) the idea of temporal dependence.

\subsection{Some Notions on the Loss Function.}

We make the following mild assumptions about the function $\ell()$:
\begin{itemize}
    \item $\ell:\mbR\rightarrow\mbR$ is a differentiable, convex, even function (so $\ell'$ is odd). 
    \item $\ell()$ has a bounded derivative ---$||\ell'||_\infty\le\tau$ for some $0\le\tau<\infty$.
    \item $\ell()$ has positive curvature on a small interval, i.e., $\ell''$ exists and satisfies $\ell''>0$ in a small neighbourhood around the origin.
\end{itemize}
If $\ell()$ is non-convex, then local/global minima may not exist and additional constraints are required (usually these constraints "cheat" by "convexifying" the problem -- the penalized procedure \eqref{stoch-reg} is now done over some local convex region in $\mbR^p$ so that local minima are guaranteed at least). A classical choice of a robust loss function is the Huber loss (and its many variants), given by
\begin{gather}\label{huber-loss}
\ell(u)=\begin{cases}\frac{u^2}{2} \quad \quad \quad \quad :|u|\le\tau \\
\tau |u|-\frac{\tau^2}{2} \quad :|u|>\tau,
\end{cases}
\end{gather}
where the "cut-off point" $\tau>0$ controls the level of robustification. Small values of $\tau$ restrict the influence of the quadratic least squares loss, hence its susceptibility to noisiness, while large values of $\tau$ do the opposite -- in the limit $\tau=\infty$, $\ell()$ is just the least squares loss. Note that the Huber loss is a convex, even function, $||\ell'||_\infty\le\tau$ and $\ell''=1$ on $(-\tau,\tau)$. In addition, $\ell''$ exists everywhere except $u=\pm\tau$.

The reason for including a weight function $\mfw():\mbR^p \rightarrow[0,\infty)$ is that, often in time series, the covariate process $\{x_t\}$ is driven by heavy shocks $\{\epsilon_t\}$, which makes $\{x_t\}$ heavy tailed also (e.g. autoregressions). A classical choice of such a weight was proposed by \citep{mallows1975some}:
\begin{gather}\label{mallows-weight}
\mfw(x)=\min\parenc{1,\frac{b}{||B_M x||}}, \quad\text{or}\quad 
\mfw(x)=\min\parenc{1,\frac{b^2}{||B_M x||^2}},
\end{gather}
where the (positive definite) shrinkage matrix $B_M$ and $b>0$ are fixed parameters. These weights essentially shrink data points $x$ for which $||x||$ is large, toward an elliptical shell defined by $B_M$; specifically $||\mfw(x)x||\le b/
\Lambda_{\min}(B_M):=b_M$. Such estimators are called \textit{bounded influence} estimators \citep{rousseeuw2005robust}. From here on, we will assume the weights $\mfw()$ bound the predictors $\{x_t\}$ in the following way: $||\mfw(x_t)x_t||\le b_M$, for $1\le t\le n$. We will call the term "$b_M$" --- Mallow's parameter (similar to Huber's truncation parameter "$\tau$"). In simulations, we take the first choice of the weight function $\mfw()$ and let the shrinkage matrix $B_M=I$ so that the Mallow's parameter $b_M=b$.

\subsection{Some Notions on the Penalty.}

As mentioned, the standard assumption behind consistently estimating a parameter living in a high dimensional space, is that, the true parameter space is actually a low-dimensional subspace or "sparse". To induce this sparsity pattern in the corresponding estimator, we use a generic penalty $\mcR()$ in \eqref{stoch-reg}. As such, some generic quantities related to $\mcR()$, will pop up in our main consistency result. To that end, we require the following definitions.

\begin{definition}\label{GaussianWidth}{Gaussian Width:}
Let G be a $p \times d$ random matrix with iid $N(0,1)$ entries. For a a set $\msT\subseteq\mbR^{p\times d}$, the Gaussian width of $\msT$ is defined as
\begin{gather}\label{gaussian-width}
w(\msT)= \underset{W\in\msT}{\sup}\tr(W^\top G) .   
\end{gather}
It measures the size of a (usually convex) set in the Euclidean space. A key challenge is to evaluate Gaussian widths of sets related to the penalty $\mcR()$.
\end{definition}

\begin{definition}\label{CompatibilityConst}{Subspace Compatibility:}
Given a set $\msC\subseteq\mbR^{p\times d}$, and a generic norm $\mcR()$ on $\mbR^{p\times d}$, the subspace norm compatibility constant is given by 
\begin{gather}\label{compatibility-constant}
    \Phi_\mcR(\msC)= \underset{W \in\msC-\{0\}}{\sup} \frac{\mcR(W)}{||W||_F}.
\end{gather}
When $d=1$, the definition reduces to the one given in \citep{negahban2012unified}. It measures the  relative price paid for switching between a generic norm $\mcR()$ and the usual Euclidean norm $||\cdot||_F$. Also, the reverse norm compatibility is given by
\begin{gather}\label{reverse-comp}
    \overline{\Phi}_\mcR (\msC)= \underset{W \in\msC-\{0\}}{\sup} \frac{||W||_F}{\mcR(W)} .
\end{gather}
\end{definition}

In particular when $\msC=\mbB_\mcR(0,1)$, we write $\overline{\Phi}_\mcR(\msC)$ as simply $\overline{\Phi}_\mcR$. For example, we note that $||v||\le\overline{\Phi}_\mcR \mcR(v)\le \overline{\Phi}_\mcR$ for all $v\in \mbB_\mcR(0,1)$. We further assume that $\overline{\Phi}_\mcR$ is bounded above by an absolute constant (this holds in most cases of interest). For more on norms that we can choose as our penalty, see Appendix \eqref{more-on-norms}.

\subsection{Some Notions on Temporal Dependence: Concentration Inequalities.}

From a technical standpoint, temporal dependence factors arise while using concentration inequalities, in course of proving consistency of our penalized estimates. Our goal is to recover \textit{optimal} consistency rates that have already been derived in the independence setting, the only difference being a dependence factor (i.e., the "price paid" for temporal dependence), which appears as a multiplicative constant in our results. Unfortunately, there is no single, master framework for temporal dependence that gives us tight concentration inequalities in every example. Hence we require different, albeit related, notions of dependence that suit specific examples. Several such notions exist, but we won't go into all of them --- for a comprehensive overview, see \citep[Supplementary, Appendix E]{basu2015regularized}. However, these different notions of dependence generally do not imply one another, and may instead, be seen as complementing each other. In this paper, we'll mainly focus on two notion of dependence --- (a) Markov Chains and (b) "mixing" (specifically $\beta$-mixing). 

\begin{definition}{Stationarity:} 
A process $\{x_t\}$ is $\textit{strictly}$ stationary if for all $m,n,r\in\mathbb{N}$, the vector $(x_m,...,x_{m+n})$ has the same distribution as $(x_{m+r},...,x_{m+n+r})$. It is $\textit{weakly}$ (or covariance) stationary if the autocorrelations $Corr(x_t,x_{t+l})$ do not depend on $t$ for all $l\in\mbZ$. For Gaussian processes, the two notions coincide. However, this fails to hold in general, e.g. heavy tailed processes. 
\end{definition}

As noted, from a technical standpoint, temporal dependence is only relevant while using concentration inequalities --- as such, it is \textit{latter} that is the main tool of interest. Since the notions of dependence require a bit of lengthy exposition, hence, to reduce clutter, this exposition itself is relegated to the Appendix, while the concentration inequalities connected to said notions, are given here.

The main tool we'll use is a concentration inequality connected to bounded functions of stationary Markov chains, satisfying some technical requirements (see Appendix \eqref{temp-dep}), given by Proposition \eqref{ConcMarkov}.

\begin{proposition}\label{ConcMarkov}

Let $f$ be a bounded function on $\mbR^p$ satisfying $|f|\le L$. Let us suppose $\{x_t\}$ is a stationary Markov chain on $\mbR^p$ satisfying a minorization \eqref{minorization} and drift \eqref{drift} condition. Further, let us suppose the random variables $T_1, T_2$ defined in the Appendix \eqref{temp-dep}, have finite exponential moments. Then, for all $t>0$, and some absolute constant $c>0$,
\begin{gather}
\mbP\parens{\biggl|\sum_{i=1}^n f(x_i)-\mbE\bigl[\sum_{i=1}^n f(x_i)\bigr] \biggr|\ge Ln\msfC_{MC}t}\le 2\exp[-cnt^2] .   
\end{gather}
\end{proposition}

This result is just Hoeffding inequality \citep[section 2.2]{vershynin2018high}, \textit{modulo} a temporal dependence factor given by the quantity $\msfC_{MC}$, which measures the exponential tails of the random lengths $T_1, T_2$ of the independent ``blocks'' defined in the Appendix \eqref{temp-dep}. 

Next, a concentration inequality that is connected to the notion of "mixing" is given in Proposition \eqref{ConcWongBounded}.

\begin{proposition}\label{ConcWongBounded}

Let $\{x_t\}$ be a (sub)-geometrically $\beta$-mixing process in $\mbR^p$ with $\beta_{mix}>0$ and geometric index $\gamma_1\in(0,2)$, and let $\gamma:=[1/\gamma_1+1/2]^{-1}$ --- note that $\gamma\in(0,1)$ (the terms $\beta_{mix}$, $\gamma_1$ are defined in Appendix \eqref{temp-dep}). Let $f$ be a bounded function on $\mbR^p$ satisfying the bound $|f|\le L$. There exists a constant $\msfC_{\mathrm{Mix}}>0$ that depends on $\beta_{mix}$ and $\gamma_1$ such that, for $n>4,t>1/n$, 
\begin{gather}
\mbP\biggl[\biggl|\sum_{i=1}^n f(x_i) -\mbE \bigl[\sum_{i=1}^n f(x_i) \bigr] \biggr|>Lnt \biggr]
\le n\exp[-\msfC_{\mathrm{Mix}}^{-2}\min\parenc{(nt)^\gamma, nt^2}] \\
\le n\exp[-\msfC_{\mathrm{Mix}}^{-2} n^\gamma\min\parenc{1,t^2}] .
\end{gather}
\end{proposition}

This follows from \citep*[Lemma 13]{wong2020lasso} and the fact that any function of a mixing process is also mixing with the same mixing coefficients. This result is very close to the classical Hanson-Wright inequality \citep[section 6.2]{vershynin2018high}, again, modulo a temporal dependence factor.

\begin{remark}

Proposition \eqref{ConcMarkov} is a clear improvement over Proposition \eqref{ConcWongBounded}, since the former holds for all $t>0$ and matches Hoeffding's inequality for bounded, independent random variables. However, to use the former, two technical conditions \eqref{drift} and \eqref{minorization} need to be verified, which can be challenging in practice. We should mention that there remains a gap in the literature, with regard to concentration inequalities for mixing processes vs Markov chains --- this is addressed in \citep{merlevede2011bernstein}. For more general, Bernstein bounds corresponding to Propositions \eqref{ConcMarkov} and \eqref{ConcWongBounded}, see \citep[Theorem 1.1, 1.2]{adamczak2015exponential} and \citep[Theorem 1]{merlevede2011bernstein} respectively. For a detailed comparison between the generalized versions, see \citep[Remark 4]{adamczak2015exponential}.

\end{remark}

\begin{remark}

In Proposition \eqref{ConcWongBounded}, the precondition $t>1/n$ is extremely mild and boils down to the sample size $n$ being bigger than a finite constant (that does not scale with dimension) in applications. Also, when the geometric index $\gamma_1 \rightarrow 2$, we have $\gamma \rightarrow 1$, and then the exponential term in the right hand side matches the Hanson-Wright bound for the independent case. However, a factor of $n$ still remains on the right hand side, which is sub-optimal. However, this is \textit{not} a big price to pay --- in fact, the exponential term usually dominates.  

\end{remark}

\section{Why is a Robust Procedure Necessary ?}

To motivate this section, let us briefly recall the main ideas behind proving the consistency of LASSO. In low dimensions when $n\gg p$, the ordinary least squares (OLS) estimator $\hat{\beta}_{OLS}$ consistently estimates $\beta^*$ in \eqref{stoch-reg}. The proof is standard --- the following is shown:

\begin{itemize}
    \item The cross-product or deviation term $X^\top\mcE/n$ in \eqref{stoch-reg-matrix-form} converges to 0 (in probability), or, in other words, concentrates around 0 "with high probability". \quad (Deviation Condition)
    \item The sample gram matrix $X^\top X/n$ is positive definite with high probability, i.e. its variational form $v^\top X^\top X v/n$ is lower bounded by an $\textit{unrestricted}$ constant $\alpha_{URE}>0$ over $\textit{all}$ $v\in\mbR^p$. \quad (Unrestricted Eigenvalue Condition)
\end{itemize}

\noindent
In high dimensions when $p\gg n$, and the usual LASSO (with least squares loss) is used, the Deviation Condition remains the same, but, since the sample gram matrix is now singular, the Unrestricted Eigenvalue condition can't hold. It turns out that only a slight modification is enough here --- the variational form $v^\top X^\top X v/n$ is now lower bounded, only over a \textit{restricted} set of vectors $v$ living in a subset of $\mbR^p$, by a constant $\alpha_{RE}>0$. This is the so-called "Restricted Eigenvalue" (RE) condition.

There is a balance to be struck between proving that \textit{both} the Deviation \textit{and} the RE conditions hold simultaneously with high probability. Roughly speaking, it is easy to show that the deviation term $X^\top\mcE/n$ (which is already centered) concentrates around 0 (its mean) w.h.p., when the predictors $\{x_t\}$ and noise $\{\epsilon_t\}$ are not \textit{too} noisy or heavy tailed (since that would worsen the concentration bound, hence the final consistency rates). By contrast, the opposite tendency holds while proving the RE condition --- heavier tails of $\{x_t\}$ implies the gram matrix, or its variational form $v^\top X^\top X v/n$ concentrates \textit{away from} 0 w.h.p. Therefore, proving both simultaneously means that the random processes can't be too light tailed or too heavy tailed. This issue is usually resolved in the classical iid literature either by assuming $\{x_t\}$ is deterministic (fixed design setup), or by assuming $\{x_t\}$ is exogenous, i.e. completely independent of the errors $\{\epsilon_t\}$. As such, the two conditions are usually verified separately \citep{negahban2012unified}. However, in the case of time series, especially in the case of endogenous processes, e.g. VAR$s$, this is not allowed. Both the predictors $\{x_t\}$ and errors $\{\epsilon_t\}$ are random (random design), but crucially the errors/noise $\{\epsilon_t\}$ drive both the predictors $\{x_t\}$ and responses $\{y_t\}$ (endogeneity). Hence, if the errors are heavy tailed (e.g. polynomial tails/only finitely many moments exist), then so are the predictors/responses. Therefore, proving both Deviation and RE w.h.p becomes challenging. 

In particular, under the LASSO with least squares loss, the Deviation condition becomes hard to verify in the presence of heavy tails. To simplify, suppose $\{x_t\}$ are iid subgaussian and the $\{\epsilon_t\}$ are iid from an $\alpha$-stable distribution \citep{loh2017statistical} with $\alpha<2$, and that $\{x_t\}$ is independent of $\{\epsilon_t\}$. Then, assuming $\beta^*$ is $s$-sparse, the least squares LASSO with the usual penalty $\mcR(\cdot)=||\cdot||_1$, fails to achieve the optimal $\mathcal{O}\left(\sqrt{\frac{s\log p}{n}}\right)$ rate. The reason is that the Deviation condition now \textit{fails} with high probability \citep[Lemma 2]{loh2017statistical}. Hence, we can't do without a robust procedure.

\section{The Main Result.}

Next, we make the following probabilistic assumptions regarding the joint process $\{(x_t,\epsilon_t)\}$ given in the model \eqref{stoch-reg}:
\begin{enumerate}
    \item\label{req_stationarity}\textbf{Stationarity.} The joint process $\{(x_t,\epsilon_t)\}$ is strictly stationary. 
    \item\label{req_noise}\textbf{Distributional Requirement of Noise Process.} For each $t$, the noise $\epsilon_t$ has a $\textit{conditionally symmetric}$ density, given the $\sigma$-field $\mfF_{i-1}$, where $\mfF_i=\sigma\{x_{i+1},x_i,x_{i-1},\dots\}$. This density is positive everywhere on $\mbR$. Assuming $\mbE|\epsilon_1|<\infty$, it follows that the noise $\{\epsilon_t\}$ forms a martingale difference sequence with respect to the filtration $\{\mfF_t\}$.
    \item\label{req_predictor} \textbf{Distributional Requirement of Predictor Process}: The predictor process $\{x_t\}$ satisfies either one of the following: 
     \begin{enumerate}
    \item\label{req_markov}{\textbf{Markov Chain Requirement.}} The predictor process $\{x_t\}$ is a Markov chain satisfying a minorization \eqref{minorization} and drift \eqref{drift} condition.
    \item\label{req_mixing}{\textbf{Mixing Requirement.}} The predictor process $\{x_t\}$ is (sub)geometrically  $\beta$-mixing with rate $\beta_{mix}$, and geometric index $\gamma_1$. We  define $\gamma:=\parens{1/\gamma_1+1/2}^{-1}$ and assume $\gamma<1$ (equivalently $\gamma_1<2$). The terms $\beta_{mix}$ and $\gamma_1$ are defined in Appendix \eqref{temp-dep}.
     \end{enumerate}
\end{enumerate}

\begin{remark}
All our examples satisfy Assumption \eqref{req_noise}. From a technical standpoint, this assumption is used to verify a first-order "deviation" condition, which in turn, gives us the size of the tuning parameter $\lambda_n$. Actually, rather than the full strength of Assumption \eqref{req_noise}, in practice, we only require that (a) each $\epsilon_t$ given $x_t$ has a conditionally symmetric density, positive everywhere on $\mbR$ and that (b) for any fixed vector $u\in\mbR^p$, the sequence of random variables 
\begin{gather}
w^2(x_t)\ell'[\epsilon_t w(x_t)]u^\top x_t,
\end{gather}
forms a martingale difference sequence with respect to the filtration $\{\mfF_t\}$. This in turn allows us to exploit Azuma's inequality (i.e. Hoeffding inequality for martingale difference sequences) in order to get a deviation bound with high probability.
\end{remark}

\begin{remark}
All our examples satisfy Assumption \eqref{req_predictor}. If our example satisfies \eqref{req_markov}, then we can leverage Proposition \eqref{ConcMarkov}, and if it satisfies \eqref{req_mixing}, then we can leverage Proposition \eqref{ConcWongBounded}, to obtain appropriate consistency rates. In principle, we do not assume $\textit{a priori}$ the existence of  second or even first moments of the noise $\{\epsilon_t\}$; nevertheless, in many examples, finiteness of these moments are used for easy-to-verify sufficient conditions to ensure (a) strict stationarity of the corresponding predictor and noise processes, and (b) the drift condition \eqref{drift}, which guarantees $\beta$-mixing. Instead, we consider the tail probability of the error $\mbP[|\epsilon|>\msfT/2]$, for $\msfT>0$ small enough so that $\alpha_T:=\underset{|u|\le \msfT}{\min}\ell''(u)>0$. For example, we may fix $\msfT=\tau$, so that $\alpha_\msfT=\alpha_\tau=1$ for the Huber loss with threshold $\tau$.
\end{remark}

We are now ready to present the theoretical consistency results of our sparse Huber estimator \eqref{robust-loss}, under the two regimes \eqref{req_markov} and \eqref{req_mixing}.
 
\begin{proposition}\label{det-errbd-robust-stochreg}{Theoretical consistency of sparse, robust regression under temporal dependence.}

\vspace{0.5em}
\noindent
\textbf{Consistency under Markovian Regime.}
Consider the penalized regression problem posited in \eqref{robust-loss}. Suppose that assumptions \eqref{req_stationarity},\eqref{req_noise} and \eqref{req_markov} hold. Let the noise process $\{\epsilon_t\}$ possess a finite second moment and let $\mbE[\epsilon_1^2]=\sigma_\epsilon^2$. Let the sample size $n$ and tuning level $\lambda_n$ satisfy %
\begin{gather}
\lambda_n=2\overline{\Phi}_\mcR b_M\tau \sqrt{\frac{w^2[\mbB_\mcR(0,1)]}{n}} ,\\
n\succsim \max\parenc{1,\msfC_{MC}^2}
\frac{\alpha_\msfT^2 b_M^2 w^2(\msC\cap\mbB_2)}{\alpha_{RE}^2} ,
\end{gather}
where the anti-concentration cone $\msC$ is given by
\begin{gather}
\msC=\msC(\beta^*)=\mathrm{cone}\{u:\mcR(u)/2 +\mcR(\beta^*)-\mcR(\hat{\beta})\ge 0\}    
\end{gather}
and the restricted eigenvalue $\alpha_{RE}$ is given by 
\begin{gather}
\alpha_{RE}= \frac{1}{2}\alpha_\msfT \Lambda_{\min} 
\parens{ \mbE\parens{\mfw^3(x_1)x_1 x_1^\top\mds{1}(\sigma_\epsilon^2\le\sqrt{\msfT/2})
\mds{1}(|\epsilon_1|\le\sqrt{\msfT/2}) } } .
\end{gather}
Here $\msfC_{MC}$ is the temporal dependence factor arising from the Markov chain $\{x_t\}$. The restricted eigenvalue $\alpha_{RE}$ may be seen as a scaled version of the usual choice: $\Lambda_{\min}\parens{\mbE[x_1 x_1^\top]}$. Recall that $b_M$ is the Mallow's parameter. We can then establish:
\begin{gather}
||\hat{\beta}-\beta^*||_2 \le\frac{3\lambda_n\Phi_\mcR(\msC)}{2\alpha_{RE}} , \quad
\mcR(\hat{\beta}-\beta^*) \le\frac{3\lambda_n\Phi_\mcR^2(\msC)}{2\alpha_{RE}} \quad\text{(Estimation error).}
\end{gather}

\vspace{0.5em}
\noindent
\textbf{Consistency under Mixing Regime.}
Suppose now, that assumptions \eqref{req_stationarity},\eqref{req_noise} and \eqref{req_mixing} hold. Then, let the tuning parameter $\lambda_n$ and sample size $n$ satisfy :
\begin{gather}
\lambda_n=2\overline{\Phi}_\mcR b_M\tau\sqrt{\frac{c_1 w^2[\mbB_\mcR(0,1)]}{n}} ,\\
n\succsim\parens{\msfC_{\mathrm{Mix}}^2 \max\parenc{1,\frac{b_M^4\alpha_\msfT^2}{\alpha_{RE}^2}} \frac{\msfT^2}{b_M^2}w^2(\mbB_2\cap\msC)}^{1/\gamma} ,\\
\text{with} \ \alpha_{RE}=\frac{\alpha_\msfT}{2}
\Lambda_{\min}\parens{\mbE(\mfw^3(x_1)x_1x_1^\top \mds{1}(|\epsilon_1|\le\msfT/2|))},
\end{gather}
where $\msfC_{\mathrm{Mix}}$ is the temporal dependence factor depending on the mixing rate $\beta_{mix}$ of the predictor process $\{x_t\}$, and the geometric index $\gamma_1$. We can then establish:
\begin{gather}
||\hat{\beta}-\beta^*||_2 \le\frac{3\lambda_n\Phi_\mcR(\msC)}{2\alpha_{RE}} , \quad
\mcR(\hat{\beta}-\beta^*) \le\frac{3\lambda_n\Phi_\mcR^2(\msC)}{2\alpha_{RE}} \quad\text{(Estimation error).}
\end{gather}
\end{proposition}

\begin{remark}

The issue with using a weight function $\mfw(\cdot)$ to effectively bound the predictors $\{x_i\}$, is that the restricted eigenvalue $\alpha_{RE}$ (which measures the ``curvature'' of the loss function) becomes smaller. In the case of robust regression, $\alpha_{RE}\propto\Lambda_{\min}\parens{\mbE(w^3(x_1)x_1x_1^\top \mds{1}(|\epsilon_1\le\msfT/2|)})\ll \Lambda_{\min}(\mbE(x_1x_1^\top))$, if $\{x_t\}$ has finite second moment. However ,the estimation error $||\hat{\beta}-\beta^*||$ scales inversely as $\alpha_{RE}$. Thus, bounding or truncating the predictors too aggressively via the weight function $\mfw(\cdot)$ is not recommended. On the other hand, instead of bounding the $x's$, the weights could be chosen judiciously so that $\mfw(x)x$ becomes Subweibull i.e. having a finite exponential moment, instead of subgaussian/bounded. In that case our results would still go through. If we assume the predictors $\{x_i\}$ are Subweibull to begin with, then nothing is lost by taking $\mfw\equiv 1$.

\end{remark}

\begin{remark}

The exponent $1/\gamma$ that appears in the finite sample bounds in Proposition \eqref{det-errbd-robust-stochreg} under the Mixing requirement \eqref{req_mixing}, satisfies $1/\gamma=1/\gamma_1+1/2>1$ (by assumption), where $\gamma_1$ is the (sub)geometric exponent of the mixing process $\{(x_t,\epsilon_t)\}$ (since robustification basically truncates the corresponding random processes, we cross back to the subgaussian world). These exponents appear due to applying the concentration inequality \eqref{ConcWongBounded}, and make the finite sample bounds sub-optimal. For example, suppose that $\{(x_t,\epsilon_t)\}$ is geometrically $\beta$-mixing, with geometric exponent $\gamma_1=1$ and the penalty $\mcR()$ is the $\ell_1$ norm. When the parameter $\beta^*$ in model \eqref{stoch-reg} is $s$-sparse, by Proposition \eqref{det-errbd-robust-stochreg} under Assumption \eqref{req_mixing}, the minimum sample size  $n$ scales as $n\succsim\max\parenc{(\log p)^2, (s\log p)^{3/2}}$ (modulo temporal dependence factors). This is somewhat disappointing, as the aim of a robust procedure is to recover $\textit{both}$ optimal consistency rates $\textit{and}$ optimal finite sample bounds that match those in the sub-gaussian case (modulo dependence factors). This is certainly true in the independent case \citep{loh2017statistical}. However, under the mixing framework, the finite sample bound takes a hit. This drawback is inherent in the concentration bound given by Proposition \eqref{ConcWongBounded} itself, and can't be overcome under the mixing framework.

\end{remark}

\begin{remark}

In Proposition \eqref{det-errbd-robust-stochreg}, the key algebraic quantities of interest are the following:
\begin{itemize}
    \item The Gaussian width of the unit norm ball : $w(\mbB_\mcR(0,1))$.
    \item The Gaussian width of the spherical cap of the tangent cone $\msC$ : $w(\msC\cap\mbB_2)$.
    \item The subspace compatibility constant $\Phi_\mcR(\msC)$. 
    \item The reverse compatibility constant : $\overline{\Phi}$. 
\end{itemize}
For example, in case of the usual elementwise sparsity, i.e. when $\mcR(\cdot)=||\cdot||_1$, and $\beta^*$ in \eqref{stoch-reg} is $s$-sparse, then $w(\mbB_\mcR(0,1))=\mcO[\sqrt{\log p}]$,  $w(\msC\cap\mbB_2)=\mcO[\sqrt{s\log p}]$, $\Phi_\mcR(\msC)=\mcO(\sqrt{s})$ and $\overline{\Phi}=\mcO(1)$. In other words, we recover the optimal consistency rate $\sqrt{\frac{s\log p}{n}}$, modulo a temporal dependence factor.

\end{remark}

\begin{remark}

Consider a VAR(1) process $z_t=B^\top z_{t-1}+\epsilon_t$ in $\mbR^p$ where our goal is to estimate the $p\times p$ transition matrix $B$ by splitting it up into $p$ regressions as given by \eqref{VAR-stoch-reg} and applying the robust procedure \eqref{robust-loss} to each column $B_j$, $1\le j\le p$ of the transition matrix $B$. Then the true parameter $\beta^*$ in \eqref{VAR-stoch-reg} is the $j^{th}$ column $B_j$, of the transition matrix $B$, where the dependence on the subscript $j$ has been temporarily suppressed. In that case, the anti-concentration cone $\msC$ in \eqref{det-errbd-robust-stochreg} actually depends on the subscript $j$. However, if we assume that all the columns $B_j$ of the transition matrix $B$ share the \textit{same} sparsity pattern --- e.g they are all elementwise sparse with at most $s$ many non-zero entries out of $p$ entries, then, the two algebraic quantities that depend on the cone $\msC$ --- namely $w(\msC\cap\mbB_2)$ and $\Phi_\mcR(\msC)$ can all be \textit{uniformly} upper bounded, independent of the subscript $j$ (see previous remark). Hence, assuming common sparsity pattern across the columns of the transition matrix, we will continue to use the same notation for the anti-concentration cone $\msC$ and suppress its dependence on the subscript $j$.

\end{remark}

\section{Examples.}

\begin{example}\label{VAR-weak-het-noise}
{VAR models generated by weakly heteroskedastic error processes.}

Following \citep*[section 4]{liebscher2005towards}, we consider the process $\{Z_t\}\in \mbR^p$ given by a vector Autoregressive Conditional Heteroskedastic (ARCH) model with 1-lag:
\begin{gather}\label{arch-VAR}
Z_t= B^\top Z_{t-1} +\Sigma(Z_{t-1})\eta_t, \quad t=1,\dots,n    
\end{gather}
where the random vectors $\{\eta_t\}$ are iid with symmetric continuous densities and positive everywhere on $\mbR^p$ (e.g., the family of continuous scale mixtures of Gaussian distribution). Suppose the noise $\{\eta_t\}$ have finite second moment, and that the noise variance is the Identity matrix. For ease of presentation, denote  $\Sigma():\mbR^p:\rightarrow \mbR^{p\times p}$ as the conditional variance that satisfies: $||\Sigma(\cdot)||_2$, $||\Sigma^{-1}(\cdot)||_2$ and $\det(\Sigma(\cdot))$ is bounded on compact sets in $\mbR^p$. 
In this setting, given $Z_{t-1}$, $Z_t$ is conditionally distributed with mean $B^\top Z_{t-1}$ and variance $\Sigma(Z_{t-1})\Sigma(Z_{t-1})^\top$. As a special case, when $\Sigma(\cdot)$ is a constant, say $\Sigma$, then these conditions are trivially satisfied provided $0<\Lambda_{\min}(\Sigma)\le \Lambda_{\max}(\Sigma)<\infty$. We also assume that the spectral radius $\rho(B)=\rho(B^\top)<1$. We examine the following two regimes:
\begin{itemize}
    \item[(i)] $||\Sigma(z)||_2=o(||z||)$ ($\textit{weak}$ conditional $\textit{heteroskedasticity}$)
    \item[(ii)] $||\Sigma(z)||_2=\mcO(||z||)$ ($\textit{strong}$ conditional $\textit{heteroskedasticity}$).
\end{itemize}
Then, assuming $\{Z_t\}$ is (strictly) stationary and $\textit{weakly}$ heteroskedastic, by \citep*[Thoerem 2]{liebscher2005towards}, $\{Z_t\}$ is geometrically $\beta$-mixing ($\gamma_1=1$, so that $\gamma=[1/\gamma_1 +1/2]^{-1}=2/3$). Hence, assuming each row of $B^\top$ is $s$-sparse and splitting the model into $p$ parallel regressions, e.g. fixing a $1\le j\le p$, the $j^{th}$ regression is then given by $y_t=x_t^\top\beta^*+\epsilon_t$, $1\le t\le n$ where $y_t=Z_{tj}$, $x_t=Z_{t-1}$, $\beta^*=B_{j:}$ and $\epsilon_t= \Sigma_j^\top(Z_{t-1})\eta_t$. From \citep{liebscher2005towards}, it is then easy to check that $\{(x_t,\epsilon_t)\}$ is a stationary, geometrically $\beta$-mixing sequence where $\epsilon_t$ is conditionally symmetric given $x_t$ (by assumption on $\eta_t$). Hence, we can carry out the $p$ regressions in parallel under the robust, mixing framework. It is important to note that, a robust framework is essential in this case, as ARCH models are well known to have polynomial tails (only finitely many moments) (cf. \citep*[Example 5]{wong2020lasso}). However, while the ``mixing'' framework applies to this case, the same drawback remains for finite sample bounds: e.g. when $\mcR$ is the $\ell_1$ penalty, ignoring dependence factors and noting that $\gamma_1=1$ and $\gamma=2/3$, the sample size $n$ scales as $\max\parenc{(\log p)^2, (s\log p)^{3/2}}$, which is sub-optimal. Finally, in practice, stationarity is not guaranteed. However, since $\{Z_t\}$ is a Markov chain in this case, we can exploit the  equivalence between ergodic Markov chains and  $\beta$-mixing (\citep*[Proposition 2, Theorem 2(i)]{liebscher2005towards}), to simulate the chain until it is close to stationarity after a sufficient ``burn-in'' period (this is standard practice for Markov Chain Monte Carlo methods, for example).

Moreover, under $\textit{weak heteroskedasticity}$, we can actually derive \textit{optimal} sample bounds in this example --- we exploit the fact that $\{Z_t\}$ is a geometrically ergodic Markov chain \citep*[Theorem 2]{liebscher2005towards} which satisfies a ``drift'' and ``minorization'' condition (as explained in section \eqref{temp-dep}). The choice of drift function used in this example is $V(z)=1+||z||$ (in the proof, the Euclidean norm $||\cdot||$ may be replaced by some equivalent vector norm). The upshot is, we can use the tighter concentration inequality \eqref{ConcMarkov} instead of \eqref{ConcWongBounded}, and get the following:

\begin{proposition}\label{det-errbd-robust-VAR}

Let the sample size $n$ and tuning level $\lambda_n$ satisfy %
\begin{gather}
\lambda_n=2\overline{\Phi}_\mcR b_M\tau \sqrt{\frac{w^2[\mbB_\mcR(0,1)]}{n}} ,\\
n\succsim \max\parenc{1,\msfC_{MC}^2}
\frac{\alpha_\msfT^2 b_M^2 w^2(\msC\cap\mbB_2)}{\alpha_{RE}^2} ,
\end{gather}
and the restricted eigenvalue $\alpha_{RE}$ be given by 
\begin{gather}
\alpha_{RE}= \frac{1}{2}\alpha_\msfT \Lambda_{\min} 
\parens{ \mbE\parens{\mfw^3(Z_0)Z_0 Z_0^\top\mds{1}(||\Sigma(Z_0)||_2\le\sqrt{\msfT/2})
\mds{1}(||\eta_0||\le\sqrt{\msfT/2}) } } .
\end{gather}
Here $\msfC_{MC}$ is the temporal dependence factor arising from the Markov chain $\{Z_t\}$. Also the restricted eigenvalue $\alpha_{RE}$ may be seen as a scaled version of the usual choice: $\Lambda_{\min}\parens{\mbE[Z_0 Z_0^\top]}$. Then, denoting the columns of B as $B_1,\dots,B_p$ we have the estimates

\begin{gather}
\underset{1\le j\le p}{\max} ||\hat{B}_j-B_j|| \le \frac{3\lambda_n\Phi_\mcR(\msC)}{2\alpha_{RE}} ,\\
\underset{1\le j\le p}{\max} \mcR[\hat{B}_j-B_j] \le \frac{3\lambda_n\Phi_\mcR^2(\msC)}{2\alpha_{RE}}
\end{gather}

with probability at least $1-2p\exp[-w^2[\mbB_\mcR(0,1)]] -2p\exp[-\frac{\msfT^2}{4b_M^2}w^2[\msC\cap\mbB_2]]$.

\end{proposition}

\begin{remark}

The crucial point here is that, by treating $\{Z_t\}$ as a Markov chain, not only do we get optimal consistency rates, but optimal sample bounds --- e.g. when the penalty $\mcR$ is the $\ell_1$ norm, the Gaussian width $w(\msC\cap\mbB_2)=\mcO(\sqrt{s\log p})$, and we see that the minimum sample size $n$ scales as $s\log p$ (ignoring other constants), which matches the independent case \citep{loh2017statistical}. 

\end{remark}

\end{example}

\begin{example}\label{AR-strong-het-noise}
{Univariate AR models with strongly heteroskedastic noise.}

We consider a univariate ARCH(p) model given by
\begin{gather}\label{univariate-ARCH}
z_t= \sum_{j=1}^p b_j z_{t-j} + \sigma(z_{t-1},z_{t-2},\dots,z_{t-p})\cdot\eta_t,
\quad t=p,p+1,p+2,\dots
\end{gather}
where $\{\eta_t\}$ is a sequence of iid random variables independent of $z_0,\dots,z_{p-1}$, having a symmetric density positive everywhere on $\mbR$ with finite second moment: $\mbE\eta^2=1$ without loss of generality. The conditional variance function is $\sigma(u)=\sqrt{d_0 + \sum_{j=1}^p d_j u_j^2}$ with $d_1,\cdots,d_p\ge 0$ and $d_0>0$. Denoting by $\tilde{Z}_t=(z_t,\dots,z_{t-p+1})^\top$, we can write the ARCH model as a Markov chain; specifically, as a VAR model in $\mbR^p$ with heteroskedastic noise terms: $\tilde{Z}_t=B^\top\tilde{Z}_{t-1} + \Sigma(\tilde{Z}_{t-1})\tilde{\eta}_t$ where
\begin{gather}
B^\top= 
\begin{bmatrix}
b_1 & b_2 & \cdots b_{p-1} & b_p \\
1   & 0   & \cdots 1       &  0  \\
0   & 1   & \cdots 0       &  0  \\
\vdots &  & \ddots &  \vdots & \vdots\\ 
0   & 0   & \cdots 1       &  0  \\
\end{bmatrix} \\
\Sigma(z)=\text{diag}[\sigma(z),0,\dots,0] 
\quad \tilde{\eta}_t=(\eta_t,0,\dots,0).
\end{gather}
Note that $||\Sigma(z)||_2=|\sigma(z)|=\mcO(||z||)$ as $||z||\rightarrow\infty$ (strong heteroskedasticity). Also, it is easy to see that $||\Sigma(\cdot)||_2= \det[\Sigma(\cdot)] =|\sigma(\cdot)|$, $||\Sigma^{-1}(\cdot)||_2 =|\sigma(\cdot)|^{-1}$ is continuous on $\mbR^p$, hence bounded on compact sets.

We can then recast the model as a stochastic regression: $y_t=x_t^\top\beta^*+\epsilon_t$, by setting $y_t:=z_t$, $x_t:=\tilde{Z}_{t-1}=(z_{t-1},\dots,z_{t-p})^\top$, $\beta^* :=(b_1,\dots,b_p)^\top$ and $\epsilon_t:=\sigma(x_t)\eta_t$. Further, assume that the spectral radius $\rho(\tilde{B})<1$, where 
$\tilde{B}=B B^\top +\text{diag}[d_1,\dots,d_p]$. Note that, under ``homoskedasticity'': $d_1,\dots,d_p=0$, this criterion reduces to the usual stability criterion: $\rho(B)<1$. Then, by \citep*[Theorem 4]{liebscher2005towards}, the Markov Chain $\{\tilde{Z}_t\}$ is Geometrically ergodic, so that $\{(x_t,\epsilon_t)\}$ is geometrically $\beta$-mixing with geometric index $\gamma_1=1$. Thus, the problem falls under our robust regression framework and Proposition \eqref{det-errbd-robust-stochreg} applies with $\gamma=2/3$. However, by Proposition \eqref{det-errbd-robust-stochreg}, the finite sample bounds suffer by an exponent term $1/\gamma=3/2$. To get rid of the exponents, we can exploit the fact that $\{\tilde{Z}_t\}$ is a Markov chain where a ``drift'' and ``minorization'' condition applies (from section \eqref{temp-dep}). In fact, the only change from example \eqref{VAR-weak-het-noise}, is the choice of the drift function; the latter can be set to $V(z)=1+||z||^2$ and then it falls under the framework posited for example \eqref{VAR-weak-het-noise}. Hence, we can use the tighter concentration inequality \eqref{ConcMarkov} instead of \eqref{ConcWongBounded}, and obtain:

\begin{proposition}\label{det-errbd-robust-AR}

For the ARCH(p) model posited in \eqref{univariate-ARCH},
let the sample size $n$ and tuning parameter $\lambda_n$ satisfy   
\begin{gather}
\lambda_n=2\overline{\Phi}_\mcR b_M\tau \sqrt{\frac{w^2[\mbB_\mcR(0,1)]}{n}} ,\\
n\succsim \max\parenc{1,\msfC_{MC}^2}
\frac{\alpha_\msfT^2 b_M^2 w^2(\msC\cap\mbB_2)}{\alpha_{RE}^2} ,
\end{gather}
and the restricted eigenvalue $\alpha_{RE}$ be given by 
\begin{gather}
\alpha_{RE}= \frac{1}{2}\alpha_\msfT \Lambda_{\min} 
\parens{ \mbE\parens{\mfw^3(\tilde{Z_p})\tilde{Z_p} \tilde{Z_p}^\top\mds{1}(||\Sigma(\tilde{Z_p})||_2\le\sqrt{\msfT/2})
\mds{1}(|\eta_p|\le\sqrt{\msfT/2}) } },
\end{gather}
wherein $\msfC_{MC}$ is the temporal dependence factor arising from the Markov chain $\{\tilde{Z}_t\}$. Then, the penalized regression coefficient $\beta^*$ satisfies:
\begin{gather}
||\hat{\beta}-\beta^*||\le \frac{3\lambda_n\Phi_\mcR(\msC)}{2\alpha_{RE}} ,\\
\mcR[\hat{\beta}-\beta^*]\le \frac{3\lambda_n\Phi_\mcR^2(\msC)}{2\alpha_{RE}} ,
\end{gather}
with probability at least $1-2\exp[-w^2[\mbB_\mcR(0,1)]] -2\exp[-\frac{\msfT^2}{4b_M^2}w^2[\msC\cap\mbB_2]]$.

\end{proposition}
Note that the restricted eigenvalue coefficient $\alpha_{RE}$ corresponds to a scaled version of what appears in penalized regression with i.i.d. errors: $\Lambda_{\min}\parens{\mbE\tilde{Z_p} \tilde{Z_p}^\top}$. 
\end{example}

\begin{example}\label{VAR-strong-het-noise}
{VAR models with strongly heteroskedastic noise}

There are a number of ways that a univariate ARCH model can be generalized to a multivariate one.

\textit{Extension 1:} Consider the VAR model $Z_t= B^\top Z_{t-1} +\Sigma(Z_{t-1})\eta_t$ as in Example 
\ref{VAR-weak-het-noise}, except now the conditional variance is defined as $\Sigma(z)= \text{diag}[(f_1+z^\top F_1 z)^{1/2}, \cdots,(f_p+z^\top F_p z)^{1/2}]$, where $F_1,\cdots,F_p$ are non-negative definite matrices and $f_1,\cdots,f_p>0$ (see \citep*[Example 4]{wu2016performance}). We also assume $\{\eta_t\}$ has finite second moment. Then, we see that $||\Sigma(z)||_2^2 = \underset{1\le j\le p}{\max} (f_j +z^\top F_j z)$, $||\Sigma^{-1}(z)||_2^2 = \underset{1\le j\le p}{\min} (f_j +z^\top F_j z)^{-1}$ , and $\det[\Sigma(z)] =\prod_{j=1}^p(f_j+z^\top F_j z)^{1/2}$ ---these are continuous on $\mbR^p$ and hence bounded on compact sets. Assume that $\rho^2(B)+\underset{1\le j\le p}{\max}\rho(F_j)<1$---note that this is stronger than the usual stability criterion: $\rho(B)<1$, although it reduces to this case under homoskedasticity (plug $F_1=\cdots=F_p=0$). We can verify that Proposition \eqref{det-errbd-robust-VAR} holds in this case also. The only change is the choice of ``drift'' function: in fact it is the same choice used for verifying Proposition \eqref{det-errbd-robust-AR} in example \eqref{AR-strong-het-noise}. Once we verify the drift condition, the rest of the proof is identical to Proposition \eqref{det-errbd-robust-VAR}.

\textit{Extension 2:} We consider a more classical generalization of the conditional variance, the so-called ``BEKK'' representation after \citep{baba1990multivariate}, given by $\Sigma(z)=[C+F^\top zz^\top F]^{1/2}$, where $C$ is a positive definite matrix. Assume that the spectral radius $\rho(B B^\top + F F^\top)<1$\footnote{It is of interest to investigate if this stability criterion could be weakened.}. We can verify that Proposition \eqref{det-errbd-robust-VAR} holds in this case as well, by verifying the drift condition (based on the same drift function from example \eqref{AR-strong-het-noise}). 
\end{example}

\begin{example}\label{VAR-threshold}
{Threshold VAR models.}

Univariate threshold autoregressive models belongs to the class of nonlinear time series models, first proposed by Tong (1978). The major features of this class of models are limit cycles, amplitude dependent frequencies, and jump phenomena \citep{tsay1989testing}. A multivariate generalization of such a model may be given as follows. Let $\mcG_1,\dots,\mcG_l$ ($l$ fixed), be a finite partition of $\mbR^p$ into $l$ disjoint regions. We consider threshold vector autoregressive model given by
\begin{gather}
z_t=\sum_{j=1}^l B_j^\top \mds{1}(z_{t-1}\in\mcG_j)z_{t-1} +\eta_t,    
\end{gather}
where $\{\eta_t\}$ is an iid sequence of random vectors with finite second moment, and a continuous symmetric density that is positive everywhere on $\mbR^p$. The parameters of interest are the matrices $B_1,\dots,B_l$. We assume the partition or thresholds $\mcG_1,\dots,\mcG_l$ are known. Let us define the map $f:\mbR^p\rightarrow \mbR^{lp}$ by $f(z)=(\mds{1}(z\in\mcG_1)z^\top:\dots: \mds{1}(z\in\mcG_l)z^\top)^\top$---then $||f(z)||^2=||z||^2$. Then writing the augmented parameter matrix $B^\top= [B_1^\top:\dots:B_l^\top]$, we can re-write the model as $z_t=B^\top f(z_{t-1})+\eta_t$. Again, we can break up this model into $p$ parallel regressions as seen in examples \eqref{VAR-weak-het-noise},\eqref{VAR-strong-het-noise}. Let us suppose $\rho(B)<1$. Then we can show the chain $\{z_t\}$ is geometrically ergodic by verifying a drift and minorization condition as before, with drift function $V(z)=1+||z||^2$. Hence the results of examples \eqref{VAR-weak-het-noise}, \eqref{VAR-strong-het-noise} follow. Of course, this can be potentially extended to heteroskedastic noise, and moreover, to more general choices of the function $f(\cdot)$, giving rise to a general class of $\textit{functional}$ vector autoregressive models.

\end{example}

\begin{example}\label{VAR-RCA} 
{Random coefficient VAR models.}

A key assumption of all standard time series models is that all parameters of the data generating model are constant (or stationary) across the observed time period. While this assumption is standard, changes of parameters over time are often plausible in psychological phenomena, especially from a $\textit{within-subject}$ perspective. For example, in the network approach to psychopathology, it is suggested that mental disorders arise from causal interactions among symptoms. These causal interactions can vary over time. Assuming there are $p$ symptoms acting as predictors, $\textit{random}$ $\textit{coefficient}$ VAR models in $\mbR^p$ capture these changes in causal interactions over time (see \citep{haslbeck2020tutorial} for more examples). Formally, a random coefficient autoregressive model \citep*[section 16.5.1]{meyn2012markov} in $\mbR^p$ is given by 
\begin{gather}
z_t=(B^\top+\Gamma_t)z_{t-1} +\eta_t ,\quad t=1,\dots,n,
\end{gather}
where $\{\Gamma_t\}$ is an iid mean zero sequence of $p\times p$ random matrices, independent of the iid mean zero noise $\{\eta_t\}$. Assume that $\mbE[\eta_0\eta_0^\top]=G$, $\mbE[\Gamma_0 \otimes\Gamma_0]=C$, where $G$ and $C$ are positive definite matrices in $\mbR^{p\times p}$ and $\mbR^{p^2\times p^2}$, respectively. The joint random vectors $\{(vec(\Gamma_t),\eta_t)\}$ have a continuous, symmetric density, positive everywhere on $\mbR^{p^2+p}$. Further, assume that the spectral radius satisfies $\rho(B^\top\otimes B^\top+C)<1$ (stability). Then, we can re-write this as a linear VAR model: $z_t=B^\top z_{t-1}+\tilde{\eta}_t$, where $\tilde{\eta}_t=\Gamma_t z_{t-1}+\eta_t$ are symmetric (dependent) random vectors with finite second moment. The results of examples \eqref{VAR-weak-het-noise}, \eqref{VAR-strong-het-noise} go through in this case also. The only non-trivial part of the verification process is the drift condition, but that is worked out in \citep*[Section 16.5.1]{meyn2012markov}. 

\end{example}

\section{Numerical Experiments.}

To back up our theory with simulations, we consider a VAR(1) process $z_t=B^\top z_{t-1}+\epsilon_t$, with $1\le t\le n$, in $\mbR^p$, where we estimate the $p\times p$ transition matrix $B=[B_1,\dots,B_p]$ by applying the robust procedure \eqref{robust-loss} to each column $B_j\in\mbR^p$ of the transition matrix $B$. We report how the estimation error $\max_{1\le j\le p}||\hat{B}_j-B_j||$ behaves in different contexts. Following \citep{loh2017statistical}, we took the Mallows shrinkage matrix to be the identity matrix $I$ and the Mallows weight $b=3$ given in \eqref{mallows-weight}, in each context.

\vspace{.1in}
\noindent
\textbf{Data Generating Process.} We chose the noise process $\{\epsilon_t\}$ in $\mbR^p$, such that, all the individual entries $\{\epsilon_{it}\}$ are iid Student's $t$ random variables, with at least 2 degrees of freedom Hence, they are automatically centered and have finite variance. We generated the transition matrix $B$ as a  adjacency matrix of an Erdos-Renyi graph and re-scaled it so that its spectral radius was $0.5$, and its sparse density (percent of non-zero edges) was $5\%$ (note that overall elementwise sparsity of the transition matrix $B$, ensures each column $B_j$ is elementwise sparse).

\vspace{.1in}
\noindent
\textbf{Estimation Procedure.} We briefly describe below, the algorithm to estimate the transition matrix $B=[B_1,\dots,B_p]$, where $B_j$ denotes the $j^{th}$ column of B. We split up problem into $p$ parallel regressions as in \eqref{intro}, and use the standard "proximal gradient descent" method, as in \citep{loh2018scale}, to estimate the parameter $B_j$ of the $j^{th}$ regression. The "proximal" or "soft-thresholding" function, defined on $\mbR^p$, is given by $S_{\alpha}(\beta)=(S^1_{\alpha}(\beta),\dots,S^p_{\alpha}(\beta))$, is defined component-wise according to
\begin{gather}
S^j_{\alpha}(\beta) = \text{sign}(\beta_j)\left(|\beta_j| - \alpha \right)_+, \qquad 1\le j\le p, \quad \alpha>0.
\end{gather}
For some step size "$\text{step}=1/\zeta$", say $\text{step}=0.9$, and some tolerance "$\text{tol}$", say $\text{tol}=10^{-4}$, and $1\le j\le p$, we build the sequence $\{B^t_j\}_{t\ge 0}$ in the following way --- we initialize $B^0_j\in\mbR^p$ by simulating its co-ordinates from an iid $\text{Uniform}(-1,1)$ distribution and then normalizing the vector. Then, for $t\ge 0$, we update $B^t_j\rightarrow B^{t+1}_j$ by the following equation: 
\begin{gather}
{B_j}^{t+1}\in\arg\min_{B_j} \left\{\mcL_n(B_j^t) +\mcL_n(B_j^t)^\top (B_j- B_j^t) + \frac{\zeta}{2} \|B_j- B_j^t\|_2^2 + \lambda\|B_j\|_1\right\} \\
=\arg\min_{B_j}\left\{\frac{1}{2} \left\|B_j- \left(B_j^t - \frac{1}{\zeta} \nabla \mcL_n(B_j^t)\right)\right\|_2^2 + \frac{\lambda}{\zeta} \|B_j\|_1\right\} \\
= S_{\lambda/\zeta}\left[B_j^t - \frac{1}{\zeta} \nabla \mcL_n(B_j^t)\right], \qquad 1\le j \le p.
\end{gather}
We stop when $||B^{t+1}_j-B^t_j||\le\text{tol}$, and take our final estimate $\Hat{B_j}=B^{t+1}_j$, for $1\le j\le p$. Hence our final estimated transition matrix is $\Hat{B}=[\Hat{B}_1,\dots,\Hat{B}_p]$.

\vspace{.1in}
\noindent
\textbf{Analyzing the Behaviour of the Error Estimate $\Hat{B}$.} Now, we report how the estimation error $\max_{1\le j\le p}||\hat{B}_j-B_j||$ behaves in different contexts. In each case, our simulations match theory, specifically, that which is given by \eqref{det-errbd-robust-stochreg} and \eqref{det-errbd-robust-VAR}.

\begin{figure}
\caption{ \label{fig:figHuberdf} \footnotesize The estimation error $\max_{1\le j\le p}||\hat{B}_j-B_j||$ plotted against various levels of noisiness (degrees of freedom). Lower degrees of freedom indicate heavier tails of the noise process. In both small VAR (figures (a) and
(c)), and Medium VAR (figures (b) and (d)), heavier robustification ($\tau=1$) gives a clear advantage over a modest one ($\tau=10$).}
\includegraphics[scale=.67]{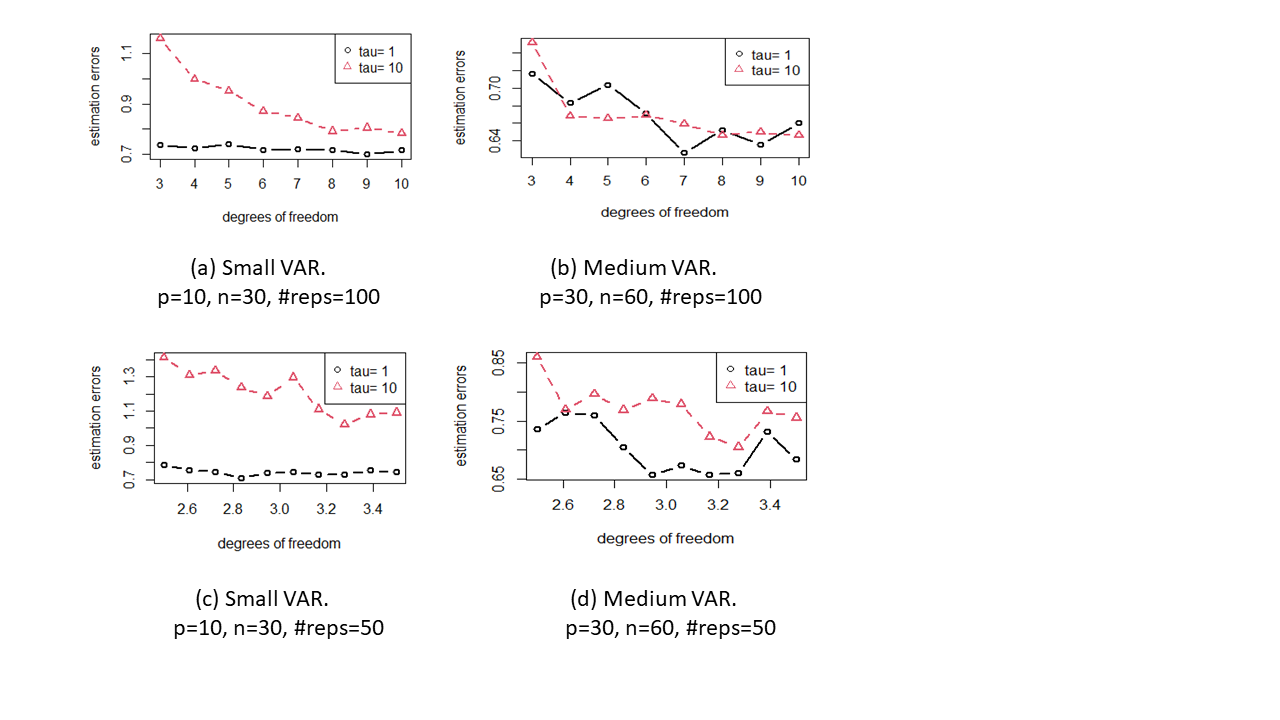}
\end{figure}

\vspace{.1in}
\noindent
\textbf{Case 1. Behaviour of $\max_{1\le j\le p}||\hat{B}_j-B_j||$ with varying degrees of freedom and fixed $n,p$.} Here we fix the dimension $p$ of the VAR(1) model (note that there are $p^2$ parameters), and sample size $n$, and see how the estimation error varies with the degrees of freedom of the $t$-noise, at two different levels of robustification --- $\tau=1$ and $\tau=10$. Figure \eqref{fig:figHuberdf} shows the results for Small VAR ($p=10, n=30$) and Medium VAR ($p=30, n=60$). When the degrees of freedom is chosen in the range $\{3,4,\dots,10\}$, the advantage of heavy robustification ($\tau=1$) becomes clear, as the tails of the noise become progressively heavier (smaller degrees of freedom). This is clearest in Figure \eqref{fig:figHuberdf}(a) --- for small VAR, the estimation error remains roughly the same at $\tau=1$ throughout, but blows up at $\tau=10$ for smaller degrees of freedom. In Figure \ref{fig:figHuberdf}(b) the advantage is much less clear --- for medium VAR, the estimation error is, on average, lower, but comparable, at level $\tau=1$, than at $\tau=10$. So, we restricted the degrees of freedom to a much smaller range $(2.5,3.5)$ (heavy noise) and found a clear advantage --- heavier robustification gives a smaller estimation error, uniformly, for both Small VAR (Figure \eqref{fig:figHuberdf}(c)), and Medium VAR (Figure \eqref{fig:figHuberdf}(d)).

\begin{figure}[p]
%
\caption{ \label{fig:figHubersample} \footnotesize The estimation error $\max_{1\le j\le p}||\hat{B}_j-B_j||$ plotted against sample size. In both small VAR (figure (a)), and Medium VAR (figure (b)), heavier robustification ($\tau=1$) gives a clear and uniform advantage over a modest one ($\tau=10$). The number of replications was chosen to be 10 to save on time.}
\includegraphics[scale=.67]{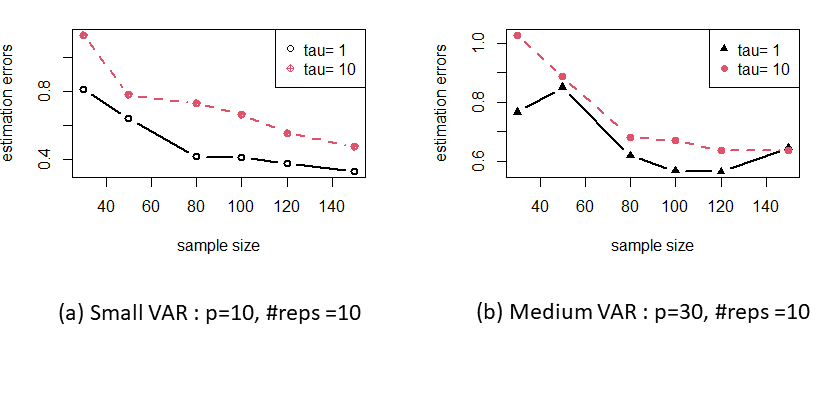}
%

%
\caption{ \label{fig:figHubersampletau} \footnotesize The estimation error $\max_{1\le j\le p}||\hat{B}_j-B_j||$ plotted against sample size. Consistency is shown for both Small and Medium VAR at each robustification level $\tau=1$ and $\tau=3$. The number of replications was chosen to be 20 to save on time.}
\includegraphics[scale=.55]{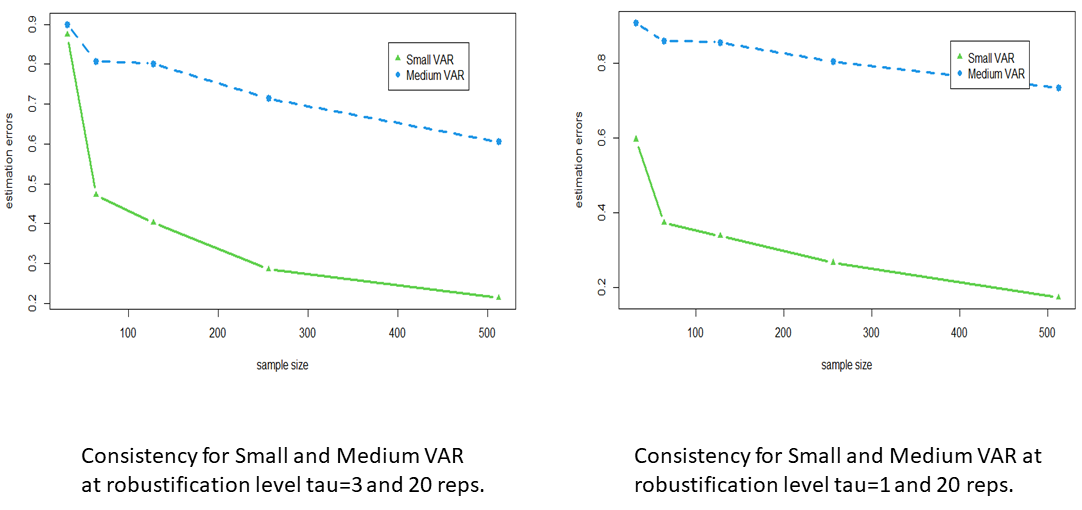}
%
    
\end{figure}

\vspace{.1in}
\noindent
\textbf{Case 2. Behaviour of $\max_{1\le j\le p}||\hat{B}_j-B_j||$ with varying sample sizes and fixed $p$ and degrees of freedom.} Here we fix the heaviness of the noise process at Student's $t$ with 3 degrees of freedom. Then we compared the behaviour of the estimation error as sample size varies, for Small VAR ($p=10$) and Medium VAR ($p=30$) --- each at 2 robustification levels $\tau=1$ and $\tau=10$. We see that in both cases (Figure \eqref{fig:figHubersample} (a) and (b)), the estimation error decreases with increasing sample size (consistency). Moreover, the estimation error is uniformly lower at a higher level of robustification ($\tau=1$) than a lower one ($\tau=10$), thereby emphasizing the importance of robustification for heavy tailed VAR.

\vspace{.1in}
\noindent
\textbf{Case 3. Behaviour of $\max_{1\le j\le p}||\hat{B}_j-B_j||$ with varying sample sizes and fixed $\tau$ and degrees of freedom.} Here we fix the heaviness of the noise process at Student's $t$ with 3 degrees of freedom. Then we compared the behaviour of the estimation error as sample size varies, for 2 robustification levels $\tau=1$ and $\tau=3$ --- each for Small VAR ($p=10$) and Medium VAR ($p=30$). We see that in both cases (Figure \eqref{fig:figHubersampletau}), the estimation error decreases with increasing sample size (consistency) for Small and Medium VAR.

\section{Discussion}

In this paper, we consider the theoretical properties of
penalized estimates in high-dimensional time series models when the data are generated from a multivariate stationary process, under structured sparsity constraints. We subsequently relaxed the assumption of Gaussian tails and replaced it with heavier tails, and thereby considered a robust approach. However, the robust approach introduces potential efficiency issues due to aggressive scaling by weight functions. Without the weight functions the finite sample bound will take a huge hit, even though we can recover optimal consistency rates. In case of heavy tails, we adopted the ``mixing'' framework to quantify temporal dependence. However this results in sub-optimal finite sample bounds, even though we recover optimal consistency rates This may be unavoidable based on the cleanest results (as far as we know) on concentration inequalities involving heavy tailed distributions under temporal dependence. It would be interesting to derive tight Bernstein or Hoeffding inequalities for mixing processes --- the one we used (Proposition \eqref{ConcWongBounded}) is the best we are aware of. However, in some cases, we showed how to re-write the time series as Markov chains and get optimal consistency rates $\textit{and}$ optimal finite sample bounds. Still, there may be interesting examples of mixing processes that cannot be written in simple Markov form, or even if that were possible, it may be difficult to verify the ``drift'' and ``minorization'' conditions required to apply the tighter concentration inequality \eqref{ConcMarkov} (e.g. multivariate GARCH models). Finally, one potential future question which we have omitted here, is (non-asymptotic) inference (based on one-step or de-sparsified estimators) for high dimensional time series under general sparsity.

\begin{appendices}



\section{More on Norms.}\label{more-on-norms}

\begin{definition}{Dual Norm:}
The dual of a norm $\mcR()$ is given by $\mcR^*(z):=\sup\{ w^\top z, \quad \mcR(w)\le 1\}$. It can be shown that $\mcR^*()$ is also a norm.
\end{definition}

\begin{definition}{Decomposable Norms:}
A norm $\mcR$ is called $\textit{decomposable}$, if for any $S\subseteq\{1,...,p\}$, for all $v\in\mbR^p$, we have $\mcR(v)=\mcR(v_S)+\mcR(v_{S^c})$.
\end{definition}

\begin{definition}{Atomic Norms:} Since many interesting norms are $\textit{not}$ decomposable, the notion of atomic norms is a useful generalization. Following \citep{chandrasekaran2012convex}, let $\mcA$ be a collection of $\textit{atoms}$ that is a compact subset of $\mbR^p$. Suppose no element $a\in\mcA$ lies in the convex hull of the other elements $conv(A-\{a\})$, i.e., the elements of $\mcA$ are the extreme points of $conv(A)$. Let $||x||_\mcA$ denote the gauge of $\mcA$, i.e. $||x||_\mcA= \inf\{t>0 :x\in t\cdot conv(\mcA)\}$.The gauge is always a convex, extended-real valued function for any set $\mcA$. By convention this function evaluates to $+\infty$ if $x$ does not lie in the affine hull of $conv(\mcA)$. We will assume without loss, that the centroid of $conv(\mcA)$ is at the origin, as this can be achieved by appropriate re-centering. With this assumption the gauge function may be recast as:
\begin{gather} \label{atomic-norm}
||x||_\mcA=\inf\{\sum_{a\in\mcA} c_a: x=\sum_{a\in\mcA}c_a a,c_a\ge 0, \forall a\in\mcA \} ,   
\end{gather}
with the sum being replaced by an integral when $A$ is uncountable. If $A$ is centrally symmetric about the origin (i.e., $a\in\mcA$ if and only if $-a\in\mcA$) we have that $||\cdot||_\mcA$ is a norm, which we call
the $\textit{atomic}$ norm induced by $\mcA$. Henceforth, we assume that $||\cdot||_\mcA$ is a norm. These norms are extremely useful because they need not be decomposable, but still cover a rich class of norms (see \citep{chandrasekaran2012convex,bhaskar2013atomic}).

\end{definition}

\section{Temporal dependence factors.}\label{temp-dep}

\subsection{Mixing.}

\begin{definition}{$\beta$-Mixing Coefficient:}
Given a probability space $(\Omega,\mcF,\mbP)$, and two sub $\sigma$-fields $\mcA,\mcB\subseteq\mcF$, the $\beta$-mixing coefficient between them \citep{bradley2005basic} is defined:
\begin{gather}
\beta(\mcA,\mcB):=\frac{1}{2}\sup \biggl\{\sum_{i=1}^I \sum_{j=1}^J |\mbP(S_i\cap T_j)-\mbP(S_i)\mbP(T_j)| ,\\
\Omega=\sqcup_{i=1}^I A_i=\sqcup_{j=1}^J B_j,\quad A_i\in\mcA,B_j\in\mcB \biggr\}
\end{gather}
where $\sqcup$ denotes a disjoint union. If two random elements $X$,$Y$ generate the $\sigma$-fields $\mcA=\sigma(X)$, $\mcB=\sigma(Y)$, then we write $\beta(X,Y)$ instead of $\beta(\mcA,\mcB)$. For a (strictly) stationary random sequence $\{(x_t,y_t)\}$, denoting $x_{-\infty:t}=\sigma\{x_j,j\le t\}$, etc. we define, for $l\in\mbZ$ 
\begin{gather}
\beta_X(l)=\beta(x_{-\infty:t},x_{t+l:\infty}),\quad
\beta_{X,Y}(l)=\beta(x_{-\infty:t},y_{t+l:\infty}),\quad 
\text{etc.}
\end{gather}
The usual practice is to define the mixing coefficients for lags $l\ge 1$.
\end{definition}

We will also use the following fact about mixing in general, repeatedly: since $\mcA'\subseteq\mcA$, $\mcB'\subseteq\mcB$ implies $\beta(\mcA',\mcB')\le\beta(\mcA,\mcB)$, for a mixing sequence $\{x_t\}$, and a $\textit{finite}$ lag $d$, the process $\{y_t\}$ defined by any measurable function $y_t:=g(x_t,x_{t-1},\dots,x_{t-d})$, is \textit{also} mixing in the same sense as $\{x_t\}$) --- with mixing coefficients bounded above by those of the original sequence $\{x_t\}$ (see e.g. \citep*[Theorem 14.1]{davidson1994stochastic}).

\begin{definition}{$\beta$-mixing Process.}
The process $\{x_t\}$ is said to be $\beta$-mixing (or $\textit{regular}$) if $\beta_X(l)\rightarrow 0$ as $l\rightarrow\infty$. Moreover a $\beta$-mixing process $\{x_t\}$ is called $\textit{sub}$-geometric if there are constants $\beta_{mix}>0$ and $0<\gamma\le 1$ such that
\begin{gather}
\beta(n)\le 2\exp\bigl\{-\beta_{mix}n^\gamma \bigr\} \quad \text{for all $n\in\mbN$ .}
\end{gather}
We call this $\gamma$ its geometric index, and $\beta_{mix}$ the mixing rate (note that $\beta_{mix}$ may depend on the true regression parameter $\beta^*$ in the model \eqref{stoch-reg}, if $\{x_t\}$ is endogenous). If $\gamma=1$, $\{x_t\}$ is said to be $\textit{geometrically}$ $\beta$-mixing. 
\end{definition}

Denoting the mixing coefficients of the \textit{joint} process $\{(x_t,y_t)\}$ by $\beta_{(X,Y)}$ (not to be confused with $\beta_{X,Y}$), it follows that $\beta_X,\beta_Y, \beta_{X,Y}\le\beta_{(X,Y)}$. \citep{wong2020lasso}, among others, prefers to work with the stronger assumption that the joint process $\{(x_t,y_t)\}$ is mixing, hence relies on coefficients $\beta_{(X,Y)}$.

\subsection{Markov Chains and Their Connection to Mixing.}

It is a common strategy to prove mixing properties of time series by first rewriting them as Markov chains, then exploiting the equivalence between ergodic Markov chains and mixing. One equivalence is the following: any stationary Markov chain is geometrically ergodic (defined below), $\textit{iff}$ it is geometrically $\beta$-mixing (see \citep*[Section 3]{bradley2005basic} for more). With that said, let $\{x_t\}$ (starting from an initial time point $t=0$), be a discrete time homogeneous Markov chain on $\mbR^p$ with Markov transition kernel $\mbP(x,A)=\mbP[x_t\in A\mid x_{t-1}=x]$ (which does not depend on $t$ by homogeneity). The ``$n$-step'' transition probability is denoted by $\mbP^n(x,A)=\mbP[x_{t+n}\in A\mid x_{t-1}=x]$. Then we have the following: 

\begin{definition}

The Markov chain $\{x_t\}$ is called ``geometrically ergodic'' if, starting from some initial point $x_0=x$, it converges to a stationary distribution $\pi$ at a geometric rate $\rho<1$. Specifically, for some function $M:\mbR^p\rightarrow [0,\infty)$,

\begin{gather}
\underset{A}{\sup} |\mbP^n(x,A)-\pi(A)|\le M(x)\rho^n, \quad n=1,2,3,\dots     
\end{gather}

\end{definition}

Geometric ergodicity is intimately connected to the notion of $\beta$-mixing--- a stationary Markov chain is geometrically ergodic $\textit{iff}$, it is geometrically $\beta$-mixing \citep*[Proposition 2]{liebscher2005towards}. Two important, classical tools used to establish geometric ergodicity are the so called ``drift'' and `minorization'' conditions. We define the concept of ``drift'' and ``minorization'' as which is now well known in literature \citep{meyn2012markov}. 

\begin{definition}\label{minorization}{Minorization condition.}

The Markov transition kernel $\mbP(x,A)$ is said to be minorized by a probability measure $\nu$ on $\mbR^p$, if there is a Borel set $\msC_{small}$ (called a ``small'' set in literature) such that, for some $\delta>0$,

\begin{gather}
\mbP(x,A)\ge\delta\nu(A) ,\quad\text{for all x}
\in\msC_{small},\quad\text{and Borel sets A}\in\mbR^p. \\
\mbP(x,\msC_{small})>0 \quad\text{for all x}\in\mbR^p.
\end{gather}

One can show that in such a situation, if the chain is stationary, with an invariant/stationary measure $\pi$, then this measure is unique and satisfies $\pi(\msC_{small})>0$ \citep{meyn2012markov}. 

\end{definition}

\begin{definition}\label{drift}{ Drift condition.}

The Markov chain  $\{x_t\}$ satisfies the ''drift''
condition, i.e. if there is a measurable function $V:\mbR^p \rightarrow [1,\infty)$, together with a constant $\lambda< 1$, such that

\begin{gather}
\mbE[V(x_t)\mid x_{t-1}=x]\le\lambda V(x) \quad\forall x\notin\msC_{small} ,\\
\underset{x\in\msC_{small}}{\sup}\mbE[V(x_t)\mid x_{t-1}=x] <\infty.
\end{gather}

The above expectation is taken with respect to the Markov transition kernel $\mbP(x,A)$. We can think of the ``drift'' function $V$ as a ``potential energy'' surface. If the ``drift'' condition holds, the chain tends to``drift'' toward states of ``lower energy'' in expectation. The function $V$ is also called the ``Lyanpunov'' function in some texts. In examples we will see, the drift function is usually taken to be $V(x)=1+||x||^s$ for some real $s>0$ (the Euclidean norm $||\cdot||$ may be replaced by an equivalent norm if necessary). Usually $s=1,2$ depending on a finite first or second moment assumption. Also, the ``small'' set $\msC_{small}$ is usually compact, so the second condition follows whenever the conditional expectation $\mbE[V(x_t)|x_{t-1}=x]$ is continuous (again, this usually holds in examples).

\end{definition}

Minorization is used to build a ``split chain'' $\{(\tilde{x}_t, R_t)\}$ on $\mbR^p\times\{0,1\}$ satisfying the following properties:

(a) $\{\tilde{x}_t\}$ is again a Markov chain with transition kernel $\mbP$ (hence for our purposes of estimating the tail probabilities we may identify $\{x_t\}$ with $\{\tilde{x}_t\}$).

(b) If we define $T_1:=\inf\{n>0:R_n = 1\}$,
$T_{i+1}:=\inf\{n>0:R_{T_1+\dots+T_{i+n}}=1\}$,
then $T_1,T_2,\dots$, form a well defined, independent random sequence.

(c) If we define $S_i:=T_1+\dots+T_i$, then the ``blocks'' $y_0 = (x_1,\dots,x_{T_1})$, $y_i:=(x_{S_{i+1}},\dots,x_{S_{i+1}})$, $i>0$, form an independent sequence. Thus, for any bounded function $f$ on $\mbR^p$, the random variables $z_i(f)
:=\sum_{j=S_{i+1}}^{S_i+1} f(x_i)$, $i\ge 1$, form an iid sequence. 

This construction is now well known in literature \citep{meyn2012markov}. For our purposes, we do not need details of this construction in its full generality. We only require that the the sequence $T_1,T_2,\dots$ (the random lengths of the independent ``blocks'') have a finite exponential moment, i.e. a moment generating function which exists finitely at least at some point. In other words, we suppose that the Subweibull norms  of $T_1$ and $T_2$ are finite, i.e. $||T_1||_{\psi_1}, ||T_2||_{\psi_1} <\msfC_{MC}<\infty$ (for a definition of Subweibull norms, see \citep[section 4]{wong2020lasso}). This assumption follows from the drift condition (see \citep*[Theorem 15.0.2]{meyn2012markov} and \citep*[section 3.5]{adamczak2008tail}). Hence, for studying time series as Markov chains, our main goal will be to verify the drift and minorization conditions. If these conditions are satisfied, we can use the construction of the ``split'' chain  and the iid ``blocks'' and get the following concentration inequality for bounded functions of geometrically ergodic Markov chains \citep*[Theorem 8]{adamczak2008tail}:

\section{Proofs of propositions.}

\subsection{Proof of Proposition \ref{det-errbd-robust-stochreg}.}

Here, we closely follow the proof technique in \citep{loh2017statistical,loh2018scale}. We first minimize over a local region $\{\beta: ||\beta-\beta^*||\le\frac{\msfT}{2b_M}\}$, i.e. we have the following (constrained) optimization problem:
\begin{gather}
\tilde{\beta}\in\underset{||\beta-\beta^*||\le\frac{\msfT}{2b_M}}{\argmin} \parenc{\mcL_n(\beta)+\lambda_n\mcR(\beta)}.
\end{gather}
We will show that this solution $\tilde{\beta}$ lies in the interior of the constraint set, hence agree with the global optima $\hat{\beta}$ of the unconstrained problem. We will prove the "algebraic" and "probabilistic" parts separately. In other words, we will assume a "deviation" and a "Restricted Eigenvalue" or "RE" condition and prove the error bounds in a deterministic fashion. Then, we will show that the "deviation" and "RE" conditions hold with high probability.
\noindent\newline
\textbf{Step 1.} Since $\tilde{\beta}$ is optimal and $\beta^*$ is feasible, the basic inequality follows:
\begin{gather}
\mcL_n(\tilde{\beta})+\lambda_n\mcR(\tilde{\beta})\le
\mcL_n(\beta^*)+\lambda_n\mcR(\beta^*)
\end{gather}
Hence, by convexity of $\mcL_n$, we get
\begin{gather}
\nabla\mcL_n(\beta^*)^\top(\tilde{\beta}-\beta^*)\le \mcL_n(\tilde{\beta})-\mcL_n(\beta^*) \le \lambda_n[\mcR(\beta^*)-\mcR(\tilde{\beta})] .
\end{gather}
Next, by Holder's inequality
\begin{gather}
0\le\mcL_n(\tilde{\beta})-\mcL_n(\beta^*) -
\nabla\mcL_n(\beta^*)^\top(\tilde{\beta}-\beta^*) \\
\le\lambda_n[\mcR(\beta^*)-\mcR(\tilde{\beta})] - \nabla\mcL_n(\beta^*)^\top(\tilde{\beta} - \beta^*) \\
\le\lambda_n[\mcR(\beta^*)-\mcR(\tilde{\beta})] + \mcR^*[\nabla_n\mcL_n(\beta^*)] \mcR(\tilde{\beta} -\beta^*).
\end{gather}
\noindent\newline
\textbf{Step 2.} Assume a deviation condition. : Here we suppose the following bound holds with high probability:
\begin{gather}
\mcR^*[\nabla\mcL_n(\beta^*)]\le\lambda_n/2 . 
\end{gather}
Then, denoting the estimate error vector by $v=\tilde{\beta}-\beta^*$, we have
\begin{gather}
0\le\mcL_n(\beta^*+v)-\mcL_n(\beta^*) -
\nabla\mcL_n(\beta^*)^\top v
\le\lambda_n[\mcR(\beta^*) -\mcR(\beta^*+v)] + \frac{\lambda_n}{2}\mcR(v).
\end{gather}
This means $v\in\msC=\msC(\beta^*)=\mathrm{cone}\{u:\mcR(u)/2 +\mcR(\beta^*)-\mcR(\beta^*+v)\ge 0\}$. Also, $\tilde{\beta}$ lies in the feasible region, so $||v||=||\tilde{\beta}-\beta^*||\le\msfT/2b_M$.
\noindent\newline
\textbf{Step 3.} Assume an RE condition. Here we suppose a (local) RE condition holds: there exists $\alpha=\alpha_{RE}>0$ such that
\begin{gather}
\mcL_n(\beta^*+u)-\mcL_n(\beta^*) -\nabla\mcL_n(\beta^*)^\top u\ge\alpha_{RE}||u||^2
\end{gather}
for every $u\in\mbB_2(\msfT/2b_M)\cap\msC$ where $\mbB_2(r')$ is the Euclidean ball of radius $r'$ around the origin. Hence, $v\in\mbB_2(\msfT/2b_M)\cap\msC$, and
\begin{gather}
\alpha_{RE}||v||^2
\le\lambda_n[\mcR(\beta^*)-\mcR(\beta^*+v)] + \frac{\lambda_n}{2}\mcR(v)
\le\frac{3\lambda_n\mcR(v)}{2} .
\end{gather}
where the last inequality is due to the triangle inequality. This finally gives
\begin{gather}
||v||=||\tilde{\beta}-\beta^*|| \le \frac{3\lambda_n \Phi_\mcR(\msC)}{2\alpha_{RE}} .
\end{gather}
Combined with feasibility of $\tilde{\beta}$, we have, thus,
\begin{gather}
||\tilde{\beta}-\beta^*|| \le\min\parenc{\frac{3\lambda_n \Phi_\mcR(\msC)}{2\alpha_{RE}} , \frac{\msfT}{2b_M}} .
\end{gather}
We will prove shortly that the tuning parameter $\lambda_n$ scales as $\sqrt{\frac{cw^2[\mbB_\mcR(0,1)]}{n}}$. Thus when $n\rightarrow\infty$ or $\lambda_n\rightarrow 0$, i.e., for $n$ large enough, or equivalently for $\lambda_n$ small enough, we see from the last inequality that the minimum of the two terms in the right hand side, just reduces to the first term:
\begin{gather}
||\tilde{\beta}-\beta^*|| \le\frac{3\lambda_n \Phi_\mcR(\msC)}{2\alpha_{RE}} .
\end{gather}
Now  we verify that the deviation and RE conditions hold with high probability, and hence, obtain lower bounds on the tuning parameter $\lambda_n$ and sample size $n$. To this end, we will use the concentration inequalities \eqref{ConcMarkov} or \eqref{ConcWongBounded}, according to whether we are following the Markovian regime \eqref{req_markov} or the Mixing regime \eqref{req_mixing}, respectively.

\noindent\newline
\textbf{Verifying the Deviation and RE condition under  the Mixing Regime.}

We verify the deviation condition using the following steps.
\noindent\newline
\textbf{Step 1.} Write the deviation term as an average of a martingale difference sequence.
\begin{gather}
\nabla\mcL_n(\beta^*)= -\frsum\ell'\parens{\epsilon_i w(x_i)} w^2(x_i) x_i.
\end{gather}
Fixing a vector $u$ so that $\mcR(u)\le 1$ (which means $||u||\le\overline{\Phi}_\mcR$), we have that the random variables
\begin{gather}
w^2(x_i)\ell'(\epsilon_i w(x_i))u^\top x_i: \quad 1\le i\le n ,
\end{gather}
are stationary, each satisfying the bound
\begin{gather}
|w^2(x_i)\ell'(\epsilon_i w(x_i))u^\top x_i| 
\le 1\cdot\tau\cdot\overline{\Phi}_\mcR||w(x_i) x_i|| \le\overline{\Phi}_\mcR b_M\tau .
\end{gather}
Moreover, by Assumption \eqref{req_noise}, they form a martingale difference sequence with respect to the filtration $\{\mfF_t\}$, where $\mfF_t=\sigma\{x_{t+1},x_t,x_{t-1},\cdots\}$, since
\begin{gather}
\mbE[\ell'[\epsilon_i w(x_i)]w^2(x_i) x_i\mid\mfF_{i-1}] \\
=w^2(x_i)x_i\mbE[\ell'[\epsilon_i w(x_i)]\mid\mfF_{i-1}]=0,
\end{gather}
using the assumptions that $\ell$ is an even function (so $\ell'$ is an odd or symmetric function, and bounded by $\tau$) and $\epsilon_i$ is conditionally symmetric given $\mfF_{i-1}$. Hence, by Azuma's inequality, for $t\ge 0$,
\begin{gather}
\mbP[|u^\top\nabla\mcL_n(\beta^*)|>t]
\le 2\exp \parenc{-\frac{2nt^2}{(\overline{\Phi}_\mcR b_M\tau)^2}}.
\end{gather}
This leads to the deviation bound, for some absolute constant $c>0$,
\begin{gather}
\mbP\parens{\mcR^*[\nabla\mcL_n(\beta^*)]
>2\overline{\Phi}_\mcR b_M\tau \sqrt{\frac{cw^2[\mbB_\mcR(0,1)]}{n}}} 
\le\exp[-w^2[\mbB_\mcR(0,1)]].
\end{gather}
\noindent\newline
\textbf{Step 2.} Choose an appropriate value of the tuning parameter $\lambda_n$. Hence, to satisfy the deviation condition, we can take  
\begin{gather}
\lambda_n=2\overline{\Phi}_\mcR b_M\tau
\sqrt{\frac{cw^2[\mbB_\mcR(0,1)]}{n}} .
\end{gather}
We verify the RE condition in the following steps. 
\noindent\newline
\textbf{Step 1.} Start with the Taylor remainder
\begin{gather}
\mcT(\beta,\beta^*):=\mcL_n(\beta)-\mcL_n(\beta^*)- \nabla\mcL_n(\beta^*)^\top(\beta-\beta^*)
\end{gather}
over the set $\{||\beta-\beta^*||\le\msfT/2b_M\}\cap\msC$. We will show that
\begin{gather}
\mcT(\beta,\beta^*)\ge\frac{\alpha_\msfT}{n} \sum_{i=1}^n  w^3(x_I)(x_i^\top(\beta-\beta^*))^2 \mds{1}(A_i)    
\end{gather}
where the sets $A_i:=\{|\epsilon_i|\le\msfT/2\}$. First, we observe that
\begin{gather}
\mcT(\beta,\beta^*)=\frsum  w(x_i) \parenc{\ell\parens{(y_i-x_i^\top\beta)w(x_i)} -\ell\parens{\epsilon_i w(x_i)} -\ell'\parens{\epsilon_i w(x_i)} w(x_i)x_i^\top(\beta^*-\beta) }.
\end{gather}
Fix an $1\le i\le n$, and let $u_{1i}=(y_i-x_i^\top\beta)w(x_i)$ and $u_{2i}=\epsilon_i w(x_i)$. Then $u_{1i}-u_{2i} = w(x_i)x_i^\top(\beta^*-\beta)$ the $i^{th}$ term of the last sum is just $w(x_i)$ times $\ell(u_{1i})-\ell(u_{2i})-\ell'(u_{2i})(u_{1i}-u_{2i})$, which is non-negative, since $\ell$ is convex. Next, $|u_{1i}|,|u_{2i}|\le\msfT$ on the event $A_i=\{|\epsilon_i| \le\msfT/2\}$.
This is because
\begin{gather}
|u_{2i}|=|\epsilon_i w(x_i)|\le |\epsilon_i|\le\frac{\msfT}{2},
\end{gather}
which means
\begin{gather}
|u_{1i}|= |(y_i-x_i^\top\beta) w(x_i)|\le |\epsilon_i w(x_i)| +|w(x_i)x_i^\top(\beta-\beta^*)| \\
\le\frac{\msfT}{2} +||w(x_i) x_i|| ||\beta-\beta^*||  \\
\le\frac{\msfT}{2} +b_M\cdot\frac{\msfT}{2b_M} \le\msfT.
\end{gather}
Further, for $|u_{1i}|, |u_{2i}|\le\msfT$, we have,
\begin{gather}
\ell(u_{1i})-\ell(u_{2i})-\ell'(u_{2i})(u_{1i} - u_{2i}) = \ell''(u_i^*)(u_{1i}-u_{2i})^2 \\
\ge \alpha_\msfT (u_{1i}-u_{2i})^2 ,
\end{gather}
where $\alpha_\msfT=\underset{|u|\le\msfT}{\min}\ell''(u)$. Thus, the claim follows:
\begin{gather}
\mcT(\beta,\beta^*)\ge\alpha_\msfT\cdot\frsum \parens{(w(x_i) x_i^\top(\beta-\beta^*)}^2 w(x_i) \mds{1}(A_i).
\end{gather}
\noindent\newline
\textbf{Step 2.} Lower bound the Taylor remainder by a quadratic form of (a scaled version of) the sample gram matrix $X^\top X/n$. Rewriting $v=\beta-\beta^*$ and $\mcT(\beta,\beta^*)=\mcT(v)$, we have, for $v\in\mbB_2(\msfT/2b_M)\cap\msC$,
\begin{gather}
\mcT(v)=\mcL_n(\beta^*+v)- \mcL_n(\beta^*) - \nabla\mcL_n(\beta^*)^\top v
\ge\alpha_\msfT\cdot\frsum \parens{x_i^\top v}^2 w^3(x_i) \mds{1}(A_i) = v^\top\Gamma v , \\
\Gamma:=\alpha_\msfT\cdot\frsum w^3(x_i) \mds{1}(A_i)x_i x_i^\top,
\end{gather}
where $\Gamma$ is a scaled version of the sample Gram matrix $X^\top X/n$. 
\noindent\newline
\textbf{Step 3.} Apply Proposition \eqref{ConcWongBounded} under Assumption \eqref{req_mixing}. Now, fix a vector $v\in\mbB_2(\msfT/2b_M)\cap\msC$. Since,
\begin{gather}
|\alpha_\msfT\parens{x_i^\top v}^2 w^3(x_i) \mds{1}(A_i) |
\le \alpha_\msfT\cdot 1\cdot ||w(x_i)x_i||^2 ||v||^2
\le\alpha_\msfT\cdot b_M^2\cdot\frac{\msfT^2}{4b_M^2} 
\le\alpha_\msfT\msfT^2/4=K_0, \quad\text{say},
\end{gather}
therefore, the term $v^\top\Gamma v$ is an average of bounded, stationary, geometrically mixing random variables (with geometric index $\gamma_1$) so that, by Proposition \eqref{ConcWongBounded}, for $n\ge 4,t>1/n$,
\begin{gather}
\mbP[|v^\top\Gamma v-\mbE(v^\top\Gamma v)|>K_0 t]
\le n\exp \parens{-\frac{1}{\msfC_{\mathrm{Mix}}^2}\min \parenc{(nt)^\gamma, nt^2}} \\
\le n\exp \parens{-\frac{1}{\msfC_{\mathrm{Mix}}^2} n^\gamma\min[1,t^2]}
\end{gather}
\noindent\newline
\textbf{Step 4.} Discretize the spherical cap $\msC\cap\mbB_2(\msfT/2b_M)$ and take union bounds.
\begin{gather}
\mbP\parens{\underset{v\in\msC, ||v||=\msfT/2b_M}{\sup} |v^\top\Gamma v-\mbE(v^\top\Gamma v)|>2K_0 t} \\
\mbP\parens{\underset{v\in\mbB_2(\msfT/2b_M)\cap\msC}{\sup} |v^\top\Gamma v-\mbE(v^\top\Gamma v)|>2K_0 t } \\
\le n\exp\parens{-\frac{1}{\msfC_{\mathrm{Mix}}^2} n^\gamma \min(1,t^2) +c' \frac{\msfT^2}{4b_M^2}w^2(\mbB_2\cap\msC)} ,   
\end{gather}
where $c'>0$ is an absolute constant (Gaussian width respects scaling).
\noindent\newline
\textbf{Step 5.} Choose an appropriate finite sample bound. We choose $n$, $t$ so that
\begin{gather}
t=\frac{\msfT^2}{16b_M^2K_0}\Lambda_{\min}[\mbE(\Gamma)]
=\frac{1}{4b_M^2\alpha_\msfT}\Lambda_{\min}[\mbE(\Gamma)] ,\\   
n\ge n_{RE}:=\parens{\msfC_{\mathrm{Mix}}^2(c'+1) \max(1,t^{-2})\frac{\msfT^2}{4b_M^2}w^2(\mbB_2\cap\msC)}^{1/\gamma}.
\end{gather}
Note that $\mbE(\Gamma)=\alpha_\msfT\mbE\parens{w^3(x_1)x_1x_1^\top \mds{1}(|\epsilon_1|\le\msfT/2)}$. We take $\alpha_{RE} =\Lambda_{\min}[\mbE(\Gamma)]/2$ and assume $\alpha_{RE}>0$. Hence, 
\begin{gather}
\mcT(v)=\mcL_n(\beta^*+v)-\mcL_n(\beta^*)- \nabla\mcL_n(\beta^*)^\top v\ge v^\top\Gamma v\ge\alpha_{RE}||v||^2 ,\quad
\forall v\in\msC\cap\mbB_2(\msfT/2b_M) ,    
\end{gather}
with probability at least $1-n\exp\parenf{\frac{\msfT^2}{4b_M^2} w^2(\mbB_2\cap\msC)}$. Also, the choice of $t$ is a constant in this case, so the precondition $t\ge 1/n$ is easily satisfied for moderate values of $n$.

Putting the pieces together, we have, for
\begin{gather}
n\ge n_{RE}, \quad
\lambda_n=2\overline{\Phi}_\mcR b_M\tau \sqrt{\frac{cw^2[\mbB_\mcR(0,1)]}{n}} ,
\end{gather}
with probability at least $1-n\exp(-w^2[\mbB_\mcR(0,1)]) -n\exp\parenf{\frac{\msfT^2}{4b_M^2} w^2(\mbB_2\cap\msC)}$,
\begin{gather}
||\tilde{\beta}-\beta^*||_2 \le \min\parenc{\frac{3\lambda_n\Phi_\mcR(\msC)}{2\alpha_{RE}}, \frac{\msfT}{2b_M} } =\frac{3\lambda_n\Phi_\mcR(\msC)}{2\alpha_{RE}},   
\end{gather}
since $\lambda_n$ scales as $\sqrt{w^2[\mbB_\mcR(0,1)]/n}$, so the rate $w^2[\mbB_\mcR(0,1)]=o(n)$ ensures the last equality. This implies that $\tilde{\beta}$ lies strictly in the interior of the local ball $\{\beta:||\beta-\beta^*||\le \frac{\msfT}{2b_M}\}$ which means that $\tilde{\beta}$ is actually a global minimizer $\hat{\beta}$ of the unconstrained problem. Furthermore, any optima of the unconstrained problem must also lie in the interior of the constraint set. This ends the proof.

\noindent\newline
\textbf{Verifying the Deviation and RE condition under Markovian Regime.}
Here, the verification is virtually identical, the only difference being that, while verifying the RE condition, when we apply the single concentration bound to the scaled gram matrix $\Gamma$, we use Proposition \eqref{ConcMarkov} under Assumption \eqref{req_markov}.

\subsection{Proof of Proposition \ref{det-errbd-robust-VAR}.}

We start with the notation for the $j^{th}$ regression given by $y_t
=x_t^\top\beta^*+\epsilon_t$, where $y_t=Z_{tj}$, $x_t=Z_{t-1}$, $\beta^*=B_{j:}^\top$, $\epsilon_t=\Sigma_j^\top(Z_{t-1})\eta_t$. We will again verify a deviation and a RE condition. However, we will use sharper concentration inequalities to get the optimal sample bounds.

For the deviation bound, we start with the gradient $\nabla\mcL_n(\beta^*)$. Fix a $u\in\mbB_\mcR(0,1)$. Then

\begin{gather}
u^\top\nabla\mcL_n(\beta^*)= 
\frsum \mfw^2(x_i)\ell'(\epsilon_i\mfw(x_i)) u^\top x_i , 
\end{gather}

where the summands $\mfw^2(x_i)\ell'(\epsilon_i\mfw(x_i)) u^\top x_i$ form a Martingale difference sequence with respect to the filtration $\mfF_i=\sigma\{x_{i+1},x_i,x_{i-1},\dots\}$ since, by conditional symmetry of $\eta_{ij}$ given $\mfF_{i-1}$, and bounded symmetry of $\ell'$,

\begin{gather}
\mbE[\mfw^2(x_i)\ell'(\epsilon_i\mfw(x_i)) u^\top x_i]\mid  
\mfF_{i-1} \\
=\mfw^2(Z_{i-1})u^\top Z_{i-1}\mbE[\ell'(\eta_{ij}\mfw(Z_{i-1}))]\mid Z_0,\dots,Z_{i-1}
=0,
\end{gather}

and 

\begin{gather}
\frac{1}{n} |\mfw^2(x_i)\ell'(\epsilon_i\mfw(x_i)) u^\top x_i|
\le \frac{1}{n}\cdot 1\cdot ||u||||\mfw(x_i)x_i||\cdot |\ell'(\epsilon_i\mfw(x_i))|
\le \frac{\overline{\Phi}_\mcR b_M\tau}{n} .
\end{gather}

Hence, by Azuma's inequality

\begin{gather}
\mbP[|u^\top\nabla\mcL_n(\beta^*)|>t]
\le 2\exp \parenc{-\frac{2nt^2}{(\overline{\Phi}_\mcR b_M\tau)^2}}.
\end{gather}

This leads to the deviation bound, for some absolute constant $c>0$,

\begin{gather}
\mbP\parens{\mcR^*[\nabla\mcL_n(\beta^*)]
>2\overline{\Phi}_\mcR b_M\tau \sqrt{\frac{cw^2[\mbB_\mcR(0,1)]}{n}}} 
\le\exp[-w^2[\mbB_\mcR(0,1)]].
\end{gather}

For the RE condition, recall that we already have a bound for a fixed $v\in\mbB_2(\frac{\msfT}{2b_M})\cap\msC$:

\begin{gather}
\mcT(v)=\mcL_n(\beta^*+v)-\mcL_n(\beta^*) -\nabla\mcL_n(\beta^*)^\top v
\ge\alpha_\msfT\cdot\frsum \parens{x_i^\top v}^2 w^3(x_i) \mds{1}(A_i) = v^\top\Gamma v , \\
\Gamma:=\alpha_\msfT\cdot\frsum w^3(x_i) \mds{1}(A_i)x_i x_i^\top,
\end{gather}

where $A_i=\{|\epsilon_i|\le\msfT/2\}$. 

Now

\begin{gather}
|\epsilon_i|=|\Sigma_j(Z_{i-1})\eta_i|
\le||\Sigma_j(Z_{i-1})||\cdot||\eta_i||
\le ||\Sigma(Z_{i-1})||_2||\eta_i|| .
\end{gather}

Hence,

\begin{gather}
||\eta_i||\le\sqrt{\msfT/2}, \quad\text{and}\quad
||\Sigma(Z_{i-1})||_2 \le\sqrt{\msfT/2}
\Rightarrow ||\Sigma(Z_{i-1})||_2||\eta_i||\le\frac{\msfT}{2} \\
\Rightarrow |\epsilon_i|=|\Sigma_j(Z_{i-1})\eta_i| 
\le ||\Sigma(Z_{i-1})||_2||\eta_i||\le \frac{\msfT}{2} \\
\Rightarrow  \mds{1}(A_i)\ge \mds{1}(B_i) \mds{1}(C_i), \\
\text{where}\quad B_i:=\{||\Sigma(Z_{i-1})||_2\le\sqrt{\msfT/2}\}, \quad C_i=\{||\eta_i||\le\sqrt{\msfT/2}\}.
\end{gather}

So we have that, for a fixed $v\in \mbB_2(\frac{\msfT}{2b_M})\cap\msC$: 

\begin{gather}
\mcT(v)\ge v^\top\Gamma v\ge v^\top\Gamma_1 v ,\\ \Gamma_1:=\alpha_\msfT\cdot\frsum \mfw^3(x_i) \mds{1}(B_i)x_i x_i^\top \mds{1}(C_i) . 
\end{gather}

Now,

\begin{gather}
v^\top[\Gamma_1-\mbE(\Gamma_1)]v
=\alpha_\msfT\cdot\frsum \mfw^3(x_i)\mds{1}(B_i)(x_i^\top v)^2 \parenc{\mds{1}(C_i) -\mbE[\mds{1}(C_i)] } \\
+\alpha_\msfT\cdot\frsum \mbE[\mds{1}(C_i)]\parenc{ \mfw^3(x_i)\mds{1}(B_i)(x_i^\top v)^2 -\mbE\parens{ \mfw^3(x_i)\mds{1}(B_i)(x_i^\top v)^2} } \\
:=I+II.
\end{gather}

Here we crucially exploit the fact that $\eta_{i}$ is independent of $x_i=Z_{i-1}$ for each $i$. To bound $I$, we note that it is a sum of terms that form a Martingale difference sequence with respect to the filtration $\mfF_i=\{x_{i+1},x_i,x_{i-1},\dots\}$. Moreover each term is bounded by $2\alpha_\msfT b_M^2\msfT^2/4b_M^2n$ =$\alpha_\msfT \msfT^2/2n$. Hence we use Azuma's inequality. To bound $II$, we note that it is also a sum of bounded functions of a Markov chain, each term satisfying the same bound $\msfT^2/2n$. So we use Proposition \eqref{ConcMarkov}. First, $\{x_t\}=\{Z_{t-1}\}$ is a geometrically ergodic Markov chain in $\mbR^p$, by \citep*[Theorem 2(i)]{liebscher2005towards}. Second, we need a drift and minorization condition for $\{Z_t\}$. 

For now, let $||B||_2<1$, which is stricter than assuming the spectral radius $\rho(B)<1$. We first verify the drift condition with drift function $V(z)=1+||z||$. By assumption, there is an $R>1$ such that, for $||z||>R>1$ (or $V(z)>R+1>2$), we have $||\Sigma(z)||_2\mbE[||\eta_0|| \le (1-||B||_2)||z||/2$. Let $\msC_{small}=\{z:V(z)\le R\}$. Then for any $z\notin \msC_{small}$, we have

\begin{align}
\mbE[V(Z_{t+1})\mid Z_t=z]
& \le 1+||B||_2||z|| +||\Sigma(z)||_2\mbE[||\eta_0||] \\
& =||B||_2 V(z) + 1-||B||_2+||\Sigma(z)||_2\mbE[||\eta_0||] \\
& \le ||B||_2 V(z)+1-||B||_2+\frac{1-||B||_2}{2}||z|| \\
& =||B||_2 V(z)+1-||B||_2+\frac{1-||B||_2}{2} (V(z)-1) \\
& =\frac{1+||B||_2}{2}V(z)+\frac{1-||B||_2}{2}  \\
& \le \frac{1+||B||_2}{2}V(z)+\frac{1-||B||_2}{2} \frac{V(z)}{2} \\
& =\underbrace{\frac{3+||B||_2}{4} }_{<1} V(z) .
\end{align}

Similarly it is easy to see that, by assumption, $\mbE[V(Z_{t+1})\mid Z_t=z]$, as a function of $z$, is bounded on the compact set $\msC_{small}$. For the minorization condition, note that it is enough to verify that the Markov transition density is minorized by a density on $\mbR^p$ (equivalently a positive, Lebesgue integrable function that can be normalized to a density) over the small set $\msC_{small}$. The Markov transition density is given by

\begin{gather}
k(t|z)= f_\eta\parens{\Sigma^{-1}(z)(t-B^\top z)}\det(\Sigma(z))^{-1} ,
\quad \forall t,z\in\mbR^p.
\end{gather}

Here $f_\eta$ stands for the joint density of the iid random vectors $\{\eta_t\}$ in $\mbR^p$. By assumption, $f_\eta$ is continuous and positive everywhere on $\mbR^p$, which means $k(\cdot|z)$ is also continuous and positive everywhere on $\mbR^p$. Fixing $t\in\mbR^p$ and taking an infimum over $z\in\msC_{small}$, we have, by assumptions on $\Sigma(\cdot)$, and compactness of $\msC_{small}$,

\begin{gather}
k(t|z)\ge \parenc{\underset{z\in\msC_{small}}{\sup}\det[\Sigma(z)] }^{-1} \parenf{\underset{z\in\msC_{small}}{\inf} f_\eta\parens{\Sigma^{-1}(z)(t-B^\top z)} } 
:=\delta\tilde{k}(t) ,
\end{gather}

where $\tilde{k}$ is a positive density on $\mbR^p$, and $\delta>0$ is a (proportionality) constant. This gives, for any Borel set $A\in\mbR^p$ and $z\in\msC_{small}$,

\begin{gather}
\mbP[z,A]=\int_{A} k(t|z) \,dt \ge \delta \int_{A} \tilde{k}(t)\, dt := \delta\nu(A),  
\end{gather}

where $\nu$  is a probability measure on $\mbR^p$ with density $\tilde{k}$. Also, $\mbP[z,\msC_{small}]>0$ for any $z\in\mbR^p$. Hence, we have shown that a drift and minorization condition holds with the assumption $||B||_2<1$. In the weaker case when the spectral radius $\rho(B)<1$, we can still guarantee the existence of a matrix norm $||\cdot||^\diamond$ on $\mbR^{p\times p}$ (depending on B), induced by a vector norm on
$\mbR^p$, such that $\rho(B) \le||B||^\diamond\le\rho(B)+ \frac{1-\rho(B)}{2} =\frac{1+\rho(B)}{2}<1$ (e.g. \citep*[p347,348]{horn2012matrix}). Finally exploiting the fact that all induced matrix norms on finite dimensional spaces are equivalent, we can replace the spectral norm $||\cdot||_2$ by $||\cdot||^\diamond$ and do the same calculations as before. 
 
We can now finally apply  Proposition \eqref{ConcMarkov} to the second term $II$ and overall, we get:

\begin{align}
\mbP[|v^\top[\Gamma_1-\mbE[\Gamma_1]]v|\ge t]
&\le \mbP[|I|\ge t/2]+\mbP[|II|\ge t/2] \\
&\le 2\exp\parens{-\frac{nt^2}{2(\alpha_\msfT\msfT^2/2)^2}} + 2\exp\parens{-\frac{cnt^2}{(\msfC_{MC}\alpha_\msfT\msfT^2)^2}}\\
&\le 2\exp\parens{-c\min\parenc{1,\frac{1}{\msfC_{MC}^2}}
\frac{nt^2}{\alpha_\msfT^2\msfT^4}} ,
\end{align}
 
where $c>0$ is an absolute constant. Thus, we have as before, denoting $K_0=\alpha_\msfT\msfT^2$,

\begin{gather}
\mbP\parens{ \underset{v\in\msC,||v||=\frac{\msfT}{2b_M}}{\sup} 
|v^\top[\Gamma_1-\mbE[\Gamma_1]]v|>2K_0 t } \\
\le 2\exp\parens{-c\min\parenc{1,\frac{1}{\msfC_{MC}^2}}
\cdot nt^2+c'\frac{\msfT^2}{4b_M^2}w^2(\mbB_2\cap\msC)} . 
\end{gather}

We choose $t$ such that

\begin{gather}
2K_0t= \frac{1}{2} \Lambda_{\min}(\mbE[\Gamma_1])||v||^2 \\
\Rightarrow t= \frac{1}{16\alpha_\msfT b_M^2} \Lambda_{\min}(\mbE[\Gamma_1]) \quad\text{since}\quad ||v||=\msfT/2b_M, K_0=\alpha_\msfT\msfT^2.
\end{gather}

With this value of $t$, we choose the sample size $n$ such that

\begin{gather}
n\ge\frac{c'+1}{c}\max\parenc{1,\msfC_{MC}^2}\cdot t^{-2}\cdot
\frac{\msfT^2}{4b_M^2} w^2(\msC\cap\mbB_2) \\
=\frac{64(c'+1)}{c}\max\parenc{1,\msfC_{MC}^2}
\frac{\alpha_\msfT^2 b_M^2 w^2(\msC\cap\mbB_2)}{\Lambda_{\min}^2(\mbE[\Gamma_1])} .
\end{gather}

Also let the restricted eigenvalue $\alpha_{RE}= \Lambda_{\min}(\mbE[\Gamma_1])/2$. Then, we have that the RE condition

\begin{gather}
\mcT(v)\ge \alpha_{RE}||v||^2 \quad\forall\,v\in\msC\cap\mbB_2(\frac{\msfT}{2b_M})    
\end{gather}

holds, with probability at least $1-2\exp\parenc{-\frac{\msfT^2}{4b_M^2}w^2(\msC\cap\mbB_2)}$ .

Thus, when the tuning parameter $\lambda_n$ and sample size $n$
satisfies

\begin{gather}
\lambda_n\asymp\overline{\Phi}_\mcR b_M\tau \sqrt{\frac{w^2[\mbB_\mcR(0,1)]}{n}} ,\\
n\succsim \max\parenc{1,\msfC_{MC}^2}
\frac{\alpha_\msfT^2 b_M^2 w^2(\msC\cap\mbB_2)}{\Lambda_{\min}^2(\mbE[\Gamma_1])} ,
\end{gather}

we have, arguing as before,

\begin{gather}
||\hat{\beta}-\beta^*||\le\frac{3\lambda_n\Phi_\mcR(\msC)}{2\alpha_{RE}} ,    
\end{gather}

with probability at least $1-2\exp[-w^2[\mbB_\mcR(0,1)]]-2\exp\parens{-\frac{\msfT^2}{4b_M^2} w^2(\msC\cap\mbB_2)}$ .

\subsection{Proof of Proposition \ref{det-errbd-robust-AR}.}

It is enough to show that the drift condition holds in this case. In fact, we will consider the more general case from example \eqref{VAR-strong-het-noise}, that is, the process $\{Z_t\}$ with conditional variance $\Sigma(z)= \text{diag}[(f_1+z^\top F_1 z)^{1/2}, \dots,(f_p+z^\top F_p z)^{1/2}]$, where $F_1,\dots,F_p$ are non-negative definite matrices and $f_1,\dots,f_p>0$. We choose the drift function $V(z)=1+||z||^2$. Then 

\begin{align}
\mbE[V(Z_{t+1})|Z_t=z] 
&\le 1+\mbE||B^\top z + \Sigma(z)\eta_0||^2  \\
&\le 1+ z^\top BB^\top z + ||\Sigma(z)||_2^2 \mbE||\eta_0||_2^2 \\
&\le 1+\rho^2(B)||z||^2 + \parenf{\underset{1\le j\le p}{\max}f_j + \underset{1\le j\le p}{\max} \rho(F_j) ||z||^2 } \mbE||\eta_0||^2 \\
&=1+\parenf{\rho^2(B)+\underset{1\le j\le p}{\max} \rho(F_j)}\parenf{V(z)-1} + \underset{1\le j\le p}{\max}f_j\mbE||\eta_0||^2 \\
&:=\lambda V(z) + K,
\end{align}

where $\lambda=\rho^2(B)+\underset{1\le j\le p}{\max} \rho(F_j) <1$ by assumption and $K=1-\lambda+ \underset{1\le j\le p}{\max}f_j\mbE||\eta_0||^2\ge 0$. Thus, defining the (compact) set $\msC_{small}=\{z: V(z)\le \frac{2K}{1-\lambda}\}$, we get that, for $z\notin\msC_{small}$,

\begin{gather}
\mbE[V(Z_{t+1})|Z_t=z] \le \frac{1+\lambda}{2}V(z),    
\end{gather}

where $0<\frac{1+\lambda}{2}<1$.

\end{appendices}

\bibliography{references}

\end{document}